\documentclass{amsart}

\usepackage{amsmath}
\usepackage{amsfonts}
\usepackage{amssymb}
\usepackage{amsthm}
\usepackage{mathrsfs}
\usepackage{verbatim}

\usepackage{etoolbox}
\makeatletter
\let\ams@starttoc\@starttoc
\makeatother
\makeatletter
\let\@starttoc\ams@starttoc
\patchcmd{\@starttoc}{\makeatletter}{\makeatletter\parskip\z@}{}{}
\makeatother

\usepackage{verbatim}

\usepackage{enumerate}
\usepackage{xcolor}

\usepackage{hyperref}
\allowdisplaybreaks

\newtheorem{theorem}{Theorem}
\newtheorem{lemma}[theorem]{Lemma}
\newtheorem{corollary}[theorem]{Corollary}

\newtheorem{proposition}[theorem]{Proposition}

\theoremstyle{definition}
\newtheorem{definition}[theorem]{Definition}

\theoremstyle{remark}
\newtheorem{remark}[theorem]{Remark}

\theoremstyle{definition}

\newcommand{\CC}{\mathbb{C}}
\newcommand{\RR}{\mathbb{R}}
\newcommand{\NN}{\mathbb{N}}
\newcommand{\ZZ}{\mathbb{Z}}

\newcommand{\tweak}{\rightharpoonup}
\newcommand{\cM}{\mathcal{M}}
\newcommand{\cD}{\mathcal{D}}
\newcommand{\bH}{\mathbf{H}}
\newcommand{\bx}{\mathbf{x}}
\newcommand{\bOh}{\mathbf{0}}
\newcommand{\be}{\mathbf{e}}

\newcommand{\rd}{\mathrm{d}}
\newcommand{\Span}{\mathrm{Span}}
\DeclareMathOperator{\supp}{supp}

\numberwithin{theorem}{section}
\numberwithin{equation}{section}

\begin{document}

\title[Ancient solutions and translators of LMCF]{Ancient solutions and translators of Lagrangian mean curvature flow}

\author{Jason D. Lotay}
\address{Mathematical Institute, University of Oxford, Oxford OX2 6GG, United Kingdom.}  \email{jason.lotay@maths.ox.ac.uk}

\author{Felix Schulze}
\address{Mathematics Institute,
  University of Warwick,
	Coventry CV4 7AL,
	United Kingdom }
\email{felix.schulze@warwick.ac.uk}
      
\author{G\'abor Sz\'ekelyhidi}
\address{Department of Mathematics, Northwestern University, Evanston, IL 60208, USA}
\email{gaborsz@northwestern.edu}

\date{\today}

\begin{abstract}
Suppose that $\cM$ is an almost calibrated, exact, ancient solution of
 Lagrangian mean curvature flow in $\CC^n$. We show that if $\cM$
has a blow-down given by the  static union of
two Lagrangian subspaces with distinct Lagrangian angles that intersect  along a
line, then $\cM$ is a translator.  In
particular in $\CC^2$, all almost calibrated, exact, ancient solutions
of  Lagrangian mean curvature flow with entropy less than 3
are special Lagrangian, a union of planes, or translators. 
\end{abstract}

\maketitle

\section{Introduction}
An important problem in complex and symplectic geometry is to find
special Lagrangian submanifolds in Calabi--Yau manifolds.
Szmoczyk~\cite{Smo96} showed that the mean
curvature flow preserves the class of Lagrangian submanifolds in
Calabi--Yau manifolds, and so one can attempt to use the flow to deform
any Lagrangian
 into a special Lagrangian. The Thomas--Yau conjecture~\cite{TY02},
motivated by mirror symmetry~\cite{Thomas01}, predicts that this is
indeed possible, assuming that the initial Lagrangian 
satisfies a certain stability condition.
More recently
Joyce~\cite{Joyce15} formulated a detailed conjectural picture,
relating singularity formation along the Lagrangian mean curvature
flow to Bridgeland stability conditions.

To motivate our main result, recall that along the mean curvature flow
of zero-Maslov Lagrangians, all tangent flows at
singularities are given by unions of   minimal  Lagrangian cones, according to
Neves~\cite{Neves.ZM}. In particular all such tangent flows
are singular, or have higher multiplicity. In order to
understand how such singularities form, it is therefore crucial to
study a general class of ancient solutions of the flow, such as Type II
blow-ups.

The simplest ancient solutions are those whose blow-down at
$-\infty$ is special Lagrangian. In this case \cite[Proposition
4.5]{LambertLotaySchulze} implies that the ancient solution itself is
special Lagrangian, and in particular static.
Our main result is the following, addressing the next simplest
situation. See Section~\ref{s:prelim} for the basic definitions. 

\begin{theorem}\label{thm:translators.intro}
Let $P_1,P_2\subset \CC^n$ be Lagrangian subspaces which intersect
along a line $\ell$ and have distinct Lagrangian angles.  Let $\cM$ be
a smoothly immersed, ancient, Lagrangian Brakke flow in $\CC^n$ with uniformly bounded area ratios. Assume further that $\cM$ is exact and zero-Maslov with uniformly bounded variation of the Lagrangian angle. For $n\geq 3$ assume in addition that $\cM$ is almost calibrated.

If $\cM$ has a blow-down at $-\infty$ given by the static flow consisting of the union $P_1\cup P_2$, then $\cM$ is a translator.
\end{theorem}

\begin{remark}\label{rem:smoothly_immersed}
For the definition of a smoothly immersed Brakke flow, see Definition \ref{def:smoothly_immersed}. Note that a mean curvature flow $F:M^n \times I \to \RR^{n+m}$ where
$I \ni t \to F(\cdot ,t)$ is a smooth family of proper immersions is an example of a smoothly immersed Brakke flow. The benefit of the notion of smoothly immersed Brakke flows is that it does not require a global parametrisation and thus the condition is preserved by local smooth convergence. In Appendix \ref{app:monotonicity} we show that the weighted monotonicity formula naturally extends to this setting, allowing weights with polynomial growth.
\end{remark}

Two possibilities for the translator in
Theorem~\ref{thm:translators.intro} are the static flows given by unions of
translates of $P_1$ and $P_2$, and the non-trivial translators constructed by
Joyce--Lee--Tsui~\cite{JLT}, which play an important role in Joyce's
conjectural picture~\cite{Joyce15}. It is an interesting question whether
there are any more possibilities.

Combining Theorem~\ref{thm:translators.intro} with the work in
\cite{LambertLotaySchulze} shows the following.  (Again, see Section~\ref{s:prelim} for the definitions.) 

\begin{corollary}
For $0<T<\infty$, let $(L_t)_{0\leq t<T}$ be a smooth, properly immersed, maximal, rational, almost calibrated Lagrangian mean curvature flow in $\CC^2$ with entropy less than $3$. Then any Type II blow-up at a singular point $(x_0, T)$ is either
special Lagrangian (and given in \cite[Theorems 1.1 and
1.3]{LambertLotaySchulze}) or a non-trivial translator. 
\end{corollary}

This result  follows from Theorem~\ref{thm:translators.intro} since for a zero-Maslov ancient solution of the flow in
$\CC^2$ with
entropy less than 3 the only possible blow-downs are a union of two
planes. If the two planes
have the same Lagrangian angle, then \cite[Theorem 1.1 or Theorem
1.3]{LambertLotaySchulze} applies. If the two planes have different
Lagrangian angle, then using \cite[Proposition 4.1]{LambertLotaySchulze}
the planes must meet along a line, and
Theorem~\ref{thm:translators.intro} applies. 
An analogous result holds also for  Lagrangian mean curvature
flow in a compact Calabi--Yau surface under a suitable rationality assumption, such as in
Fukaya~\cite[Definition 2.2]{Fukaya03}, at singularities with
density less than 3.

 In singularity analysis it is
important to consider arbitrary blow-up limits of the flow, not just
those that are smooth. Theorem~\ref{thm:translators-extension} provides an
extension of Theorem~\ref{thm:translators.intro}  in
the case $n=2$ to Brakke flows obtained as blow-up limits of smooth flows.

Recently there has been great progress in classifying ancient
solutions of geometric flows such as Ricci flow and   mean
curvature flow (see e.g.~\cite{Brendle20, DHS10, BC19, CHH18,
  CHHW19}). A crucial new difficulty in our
work is that the blow-down $P_1\cup P_2$ is singular along the
line of intersection $\ell$. As a result, an approach
based on the analysis of the linearized operator on the blow-down
faces substantial difficulties. An earlier result characterizing
translators among eternal solutions to the mean curvature flow of hypersurfaces is due to
Hamilton~\cite{Hamilton95}, relying on a differential Harnack
estimate. It is not known if this approach can be extended to higher
codimension flows. Our approach is
completely different and relies on additional structure present in
the Lagrangian setting. In particular, translators are characterized by
the condition that one of the coordinate functions, $w$, is a
linear combination of 1 and the Lagrangian angle $\theta$, i.e.~$w=a +
b\theta$, see Proposition~\ref{prop:translator.height}. The function
$w$ is an ancient solution of the heat equation along $\cM$, and the
basic idea of the proof is to obtain information about $w$ through
solutions of the heat equation on the possible blow-downs. This is related to work of
Colding--Minicozzi~\cite{CMcomplexity} on  mean curvature flow, and
also to earlier works on harmonic
functions~\cite{CMHarmonic,DingHarmonic}
 and holomorphic functions~\cite{DS2}. 

To illustrate the basic ideas, suppose that $n=2$ and let $w$ be
a coordinate function vanishing on $P_1\cup P_2$. Define $z$ so that
on $\CC^2$ we have $\nabla w = J\nabla z$, i.e.~the line $\ell$ is
parallel to $\nabla z$. Let $x$ be a coordinate function vanishing on
$P_1$, but not on $P_2$ and similarly $y$ a coordinate vanishing on $P_2$ but not
on $P_1$. The ancient solutions of the
heat equation on $P_1\cup P_2$ with at most linear growth,
allowing a different smooth solution
on each plane, are spanned by $1, \theta, x, y, z, z\theta $ (see
Lemma~\ref{lem:P1P2soln}).
Here $\theta$ is simply a different constant on each
plane. At the same time, $1, \theta$ and the coordinate functions
$x,y,z,w$ are ancient
solutions of the heat equation on our ancient flow $\cM$.

The first
main step of the proof is to show that either $\cM$ is a translator,
or along a suitable sequence of
scales $t\to -\infty$, the normalized projection of $w$ orthogonal to $x,y,z$
converges to $z \theta $ on the blow-down $P_1\cup P_2$ (see
Proposition~\ref{prop:2cases}). The main technical difficulty at this stage
is that the singular set given by the line $\ell$ has codimension one
in $P_1\cup P_2$, and we need to exploit that the angle $\theta$ takes
on different values on $P_1$ and $P_2$ in order to pass solutions of
the heat equation along $\cM$ to solutions on $P_1\cup P_2$ in the
limit. This is the content of Proposition~\ref{prop:limit.heat}.

The proof of Theorem~\ref{thm:translators.intro} is completed using
Proposition~\ref{prop:linking}, based on 
the idea that if along the flow $w$ behaves like $z\theta$ at some scale,
then the flow must break into two pieces, which roughly
look like the two planes $P_1, P_2$ rotated
in such a way that their intersections with the unit spheres are
linked.  Here the fact that $\theta$ is a different constant on each plane $P_1,P_2$ is crucial.  This linking behaviour is   used to show that the flow must have a point of density
two, but the monotonicity formula  then implies that the flow is a
static union of planes.

\subsection*{Acknowledgements.}
This project grew out of discussions at the AIM workshop ``Stability in
mirror symmetry'' in December 2020. 
JDL and FS were partially supported by a Leverhulme Trust Research
Project Grant RPG-2016-174. GSz was supported in part by NSF grant DMS-1906216. We are grateful to the referees for many helpful suggestions which improved the exposition of the paper.

\section{Preliminaries}\label{s:prelim}

In this section we introduce various key definitions and notation that we shall require throughout the article.  In particular, we introduce the  set-up for our study.

\subsection{Lagrangians in \texorpdfstring{$\CC^n$}{Cn}}  We first
recall some basic definitions concerning Lagrangian submanifolds in
$\CC^n$. 

\begin{definition}
\label{dfn:ZM}
An oriented Lagrangian $L$ in $\CC^n$ is \emph{zero-Maslov} if there
exists a function $\theta$ on $L$ (called the \emph{Lagrangian angle})
so that
\[ \Omega|_L = e^{i\theta}\, dVol_L, \]
where $\Omega = dz_1\wedge\ldots \wedge dz_n$ is the standard
holomorphic volume form on $\CC^n$ and $dVol_L$ is the Riemannian
volume form of $L$. We then have $\bH=J\nabla\theta$, 
where $\bH$ is the mean curvature vector of $L$ and $J$ is the complex
structure on $\CC^n$.  We further say that $L$ is \emph{almost
  calibrated} if $\theta$ can be chosen so that  
\begin{equation*}
\sup\theta-\inf\theta\leq \pi-\epsilon.
\end{equation*}
for some $\epsilon>0$.
\end{definition}

\begin{definition}
\label{dfn:exact}
An oriented Lagrangian $L$ in $\CC^n$ is \emph{exact} if there exists a function $\beta$ on $L$ so that
\begin{equation*}
J\bx^{\perp}=\nabla\beta,
\end{equation*}
where $\bx^{\perp}$ is the normal projection of the position vector $\bx\in\CC^n$.  Equivalently,
\[
\rd\beta=\lambda|_L,
\]
where $\lambda$ is the Liouville form on $\CC^n$, which is a $1$-form on
$\CC^n$ so that $\frac{1}{2}\lambda$ is a primitive for the K\"ahler
form $\omega$ on $\CC^n$.
The Lagrangian $L$ is \emph{rational}  if the set $\lambda(H_1(L,
\mathbb{Z}))$ is discrete in $\mathbb{R}$.   An exact Lagrangian is clearly rational.
\end{definition}

\subsection{Spacetime track}
Throughout we consider a smooth, zero-Maslov, exact, ancient solution to Lagrangian mean curvature flow (LMCF)
$$(-\infty,0)\ni t\mapsto L_t\subset \CC^n$$
which evolves with normal speed given by
 $\bH$, with uniformly bounded variation of the Lagrangian angle. We assume that
 $L_t$ has uniformly
 bounded area ratios, i.e.~there exists $C>0$ such that
 \[ \sup_{x,t}
   \mathcal{H}^n(L_t\cap B(\mathbf{x},r)) \leq C r^n \text{ for all } r>0,\]
 where
 $B(\mathbf{x},r)$ is the Euclidean ball of radius $r$ about $\mathbf{x}\in\CC^n$. 
 We call
 $$\cM:=\{L_t\times \{t\}\, |\, t \in (-\infty,0)\} \subset  \CC^n\times \RR$$
 the spacetime track of the flow, and write $\cM(t)= L_t\, .$
 
 \begin{remark}
For $n>2$ we will additionally need to assume that the flow $\cM$ is almost calibrated, so that one can apply the structure theory in \cite {LambertLotaySchulze} and \cite{Neves:survey}.
 \end{remark}

 Since our focus is on planes arising as blow-ups or blow-downs, it is useful to consider them as trivial static flows as follows.
\begin{definition}
\label{dfn:planes.static}
For a pair of  $n$-dimensional planes $P_1,P_2 \subset  \CC^n$, we let $\cM_{P_1\cup P_2}$ denote the static flow corresponding to $P_1\cup P_2$.
\end{definition}

\subsection{Rescalings}

It will  be useful to perform parabolic rescalings of our flows, so we shall introduce the following notation.
\begin{definition}
\label{dfn:parabolic}
For $\lambda>0$ we shall denote the \emph{parabolic rescaling}
$$\cD_{\lambda}:\CC^n\times \RR \to \CC^n\times \RR,\quad (\bx, t) \mapsto (\lambda \bx, \lambda^2t)\, .$$
Note that for a (Lagrangian) mean curvature flow $\cM$,  we have that
$\cD_{\lambda}\cM$
is again a (Lagrangian) mean curvature flow.
\end{definition}

It turns out to be helpful to consider a further rescaling, which turns self-similarly shrinking solutions into static points of the flow.
\begin{definition}
\label{dfn:rescaled.flow}
The \emph{rescaled flow} is 
$$ \RR \ni \tau \mapsto L_\tau := e^\frac{\tau}{2} \cM(-e^{-\tau})= e^\frac{\tau}{2} L_{-e^{-\tau}} $$
which evolves with normal speed
\begin{equation}\label{eq:rescaled.flow}
 \bH + \frac{\bx^\perp}{2}\, .
\end{equation}
\end{definition}

We recall Huisken's monotonicity formula~\cite{Huisken}:
\begin{equation}\label{eq:Huisken}
\begin{split}
 \frac{d}{dt} \int_{L_t} f \rho_{\bx_0, t_0}\, d\mathcal{H}^n &= 
    \int_{L_t} (\partial_t f - \Delta f) \rho_{\bx_0, t_0}\,
    d\mathcal{H}^n\\
    &\qquad  - \int_{L_t} f\,  \bigg|\, \bH - \frac{(\bx-\bx_0)^\perp}{2 (t-t_0)}\bigg|^2 \rho_{\bx_0, t_0}\,
    d\mathcal{H}^n\, ,  
    \end{split}
\end{equation}
for $t<t_0$, where $f$ is a smooth function on $L_t$ with (uniformly) compact support, and
\[ \rho_{\bx_0, t_0}(\bx,t) = (4\pi(t_0-t))^{-n/2} \exp\left( -\frac{
      |\bx-\bx_0|^2}{4(t_0-t)}\right) \]
is the backwards heat kernel (centred at $(\bx_0,t_0)$). For an extension to the non-compact setting (suitable for the current set-up) and functions $f$ with at most polynomial growth see Proposition \ref{prop:monotonicity}.
The density of a point $(\bx_0, t_0)$
along the flow $L_t$ is
defined to be
\[ \Theta(\bx_0, t_0) = \lim_{t \nearrow t_0} \int_{L_t} \rho_{\bx_0,
    t_0}\, d\mathcal{H}^n. \]
Recall also the  \emph{entropy} $\mu(L)$ defined by Colding--Minicozzi
\cite{ColdingMinicozzi:generic}: 
 \[
  \mu(L)=\sup_{\mathbf{x}_0\in\CC^n,\,r>0}\frac{1}{(4\pi
    r)^{n/2}}\int_L
  e^{-\frac{|\mathbf{x}-\mathbf{x}_0|^2}{4r}}\rd\mathcal{H}^n. 
 \]
By virtue of Huisken's monotonicity formula, $t\mapsto \mu(L_t)$ is non-increasing along any $n$-dimensional mean curvature flow in $\CC^n$.

\subsection{Set-up}\label{ss:setup}  We now describe the main set-up that we shall have throughout the article.  In particular, this will be useful to fix notation.

 We consider two oriented Lagrangian planes $P_1,P_2 \subset  \RR^{2n} = \CC^n$   which intersect along an oriented (real) line $\ell$ through $0$. Suppose further that $P_1,P_2$ have \emph{distinct} Lagrangian angles, which we denote by $\overline{\theta}_1,\overline{\theta}_2$.  (Note that this must be the case if $n=2$.)  Changing the Lagrangian angles by a fixed constant, we can in addition assume that $\overline{\theta}_1= - \overline{\theta}_2$. 
 
 We assume that the unit vector in the direction  of $\ell$ is given by $\be_z$, corresponding to the (real) coordinate $z$. We let $\be_w = J \be_z$, corresponding to the real coordinate $w$, noting that $\be_w$ is necessarily orthogonal to $P_1\cup P_2$.  We will think of $w$ as the ``\emph{height}'' function, since $w=0$ on $P_1\cup P_2$. We choose  coordinates $x_1,\ldots,
x_{2n-2}$ such that $x_1,\ldots, x_{n-1}, w$ vanish along $P_2$ and
$x_n,\ldots, x_{2n-2}, w$ vanish along $P_1$.

Our key assumption is that our ancient solution $\cM$ to LMCF has a blow-down at $-\infty$ given by $P_1 \cup P_2$, i.e.
\begin{equation}\label{eq:blow-down} \cD_{\lambda_i}(\cM)\cap\{t<0\} \tweak \cM_{P_1\cup P_2} \cap \{t<0\}
\end{equation}
for $\lambda_i \searrow 0$, where the convergence is in the sense of Brakke flows. We note that this is equivalent to the assumption that the sequence of smooth flows
$$(-\infty,0)\ni t\mapsto L_t^i=\lambda_iL_{\lambda_i^{-2}t} $$
for $t<0$ converges weakly to the (immersed) static flow $(-\infty,0)\ni t\mapsto P_1\cup P_2$.

At this point we make the following observation about blow-downs of $\cM$.

\begin{proposition}\label{prop:blowdown.line} Let $\cM$ be an ancient solution to Lagrangian mean curvature flow as above satisfying \eqref{eq:blow-down}, so has blow-down at $-\infty$ given by $P_1\cup P_2$ and for $n\geq 3$ is almost calibrated. 
 Then all blow-downs at $-\infty$ of $\cM$ are unions of two multiplicity one planes which meet along a subspace $L$ of dimension $m \in \{1, \ldots, n-1\}$ and have Lagrangian angles $\overline{\theta}_1$, $\overline{\theta}_2$.
\end{proposition}

\begin{proof}
We have the set of angles $\{\overline{\theta}_1, \overline{\theta}_2\}$ for two multiplicity one planes in one blow-down and the set of angles is the same for any blow-down by \cite[Theorem 3.1]{LambertLotaySchulze}.  Therefore, since $\overline{\theta}_1\neq\overline{\theta}_2$, and by \eqref{eq:blow-down} the (Gaussian) density at $-\infty$ is two, any blow-down has to consist of two distinct unit multiplicity planes with the given angles. Again using $\overline{\theta}_1\neq \overline{\theta}_2$, the case that a blow-down is a pair of transverse planes is ruled out since in this case \cite[Proposition 4.1]{LambertLotaySchulze} would force the ancient solution $\cM$ to be the static flow consisting of the two transverse planes, which contradicts one blow-down being two planes meeting along a  line.
\end{proof}

\subsection{Translators}
There is a smooth, connected,  zero-Maslov (in fact, almost calibrated), exact, ancient (in fact, eternal) solution  to Lagrangian mean curvature flow in $\CC^n$ whose blow-down at $-\infty$ is $P_1\cup P_2$, which is a \emph{translator} constructed by Joyce--Lee--Tsui \cite{JLT}.  We recall the definition of a translator as follows.

\begin{definition}
 A \emph{translator} (in the $\be_z$ direction) is a solution of LMCF satisfying 
  \begin{equation}\label{eq:translator}
  \bH=\kappa\be_z^{\perp}
  \end{equation}
  for some $\kappa\neq 0$ at some (and therefore any) time. (In fact, we can rescale the translator and change its orientation   so that any $\kappa\neq 0$ can be realised).
\end{definition}

Suppose we have a zero-Maslov translator satisfying \eqref{eq:translator}. Since $\bH=J\nabla\theta$, where $\theta$ is the Lagrangian angle, $\be_z=-J\be_w$ and thus $\be_z^{\perp}=-J\be_w^{\rm T}$, we deduce that
\begin{equation}\label{eq:translator.angle}
 \theta+\kappa w=c
 \end{equation}
for some $\kappa\neq 0$ and constant $c$ (on each component of the translator).  

\begin{remark}
It is worth noting that any Lagrangian plane $P$ so that $\be_z$ is tangent to $P$ will give a trivial example of a Lagrangian translator, since it will satisfy \eqref{eq:translator}.  
\end{remark}

Using \eqref{eq:translator.angle} we deduce the following  result.

\begin{proposition}
\label{prop:translator.height}
If the flow $\cM$ satisfies \eqref{eq:blow-down}, so has blow-down at $-\infty$ given by $P_1\cup P_2$ and for $n\geq 3$ is almost calibrated, then it is a translator in the $\be_z$ direction if and only if on each component of $\cM$ the height satisfies
\begin{equation}\label{eq:height.translator}
w=a+b\theta
\end{equation}
for some constants $a$ and $b$. 
\end{proposition}

\begin{proof} If $\cM$ is a translator in the $\be_z$ direction then \eqref{eq:height.translator} is satisfied by \eqref{eq:translator.angle}.  

We now suppose that \eqref{eq:height.translator} is satisfied on $\cM$.
If $b\neq 0$ we deduce that $\cM$ is a translator by differentiating \eqref{eq:height.translator} along $\cM(t)$ for each $t$, which yields $\bH=b^{-1}\be_z^{\perp}$.  

If $b=0$ then $w$ is constant on each component of $\cM$ and so $\be_z$ is tangent to $\cM$ as the flow is Lagrangian. Hence, $\cM$ splits as   $\cM'\times \RR$, where  $\cM'$ is an ancient solution to Lagrangian mean curvature flow in $\CC^{n-1}$.

If $n\geq 3$, then $\cM'$ is almost calibrated and by \eqref{eq:blow-down} has a blow-down at $-\infty$ given by $P'_1\cup P'_2$ (where $P_i= P_i'\times \ell$). Note that $P'_1, P'_2$ are transverse, but have different Lagrangian angle. Then \cite[Proposition 4.1]{LambertLotaySchulze} implies that $\cM' = \cM'_{P'_1\cup P'_2}$ and thus $\cM = \cM_{P_1\cup P_2}$.

 If $n=2$, then  $\cM'$ is an ancient solution $\gamma$ to curve shortening flow in $\RR^2$, which has a blow-down at $-\infty$  which is a pair of non-parallel lines.   We now show that $\gamma$ must in fact be the asymptotic lines.

\begin{lemma}
\label{lem:static.curve}
Let  $\gamma = (\gamma(t))_{-\infty<t<T}$ be an ancient smooth curve shortening flow in $\RR^2$. Assume that a blow-down $\mathcal{D}_{\lambda_i}(\gamma)$ (for $\lambda_i \searrow 0$) is either a pair of unit density static lines  $\ell_1 \cup \ell_2$ meeting at one point or a single unit density line. Then the  flow $\gamma$ is the static line(s).
\end{lemma}

\begin{proof}  The case where the blow-down is a single unit density line follows from the monotonicity formula, so we only consider the case of a pair of transverse lines in the blow-down.

  If $\gamma$ were almost calibrated, then the classification of almost calibrated ancient solutions to Lagrangian mean curvature flow in \cite[Proposition 4.1]{LambertLotaySchulze}  implies that $\gamma$ must be the lines, since they have distinct angles.

 Let  $\gamma_i:=\mathcal{D}_{\lambda_i}(\gamma)$.  Using the
 right-hand side of the integrated monotonicity formula
 (\eqref{eq:Huisken} with $f =1$) and Fatou's lemma we can pick a time
 $t < 0$ (and a subsequence in $i$) such that on $\gamma_i(t)$ the
 curvature of the curve is locally uniformly bounded in $L^2$.  This
 implies locally uniform convergence in $C^{1,\alpha}$  from which it
 follows that, for $i$ and $R$ sufficiently large, $\gamma_i(t)\cap
 B_R(0)$ is given by the union of two small $C^{1,\alpha}$ graphs over
 $(\ell_1 \cup \ell_2)(t)$. Since the flow is smooth (and using the
 pseudolocality result \cite[Theorem 1.5]{IlmanenNevesSchulze}), this description of the flow has to persist for a short time. We deduce that the flow has a point with Gaussian density two (where the two graphs intersect) and thus the flow is backwards self-similar around that point. Since we have assumed that one blow-down is $\ell_1 \cup \ell_2$, the result follows.  
\end{proof}

By Lemma \ref{lem:static.curve} we deduce that each component of $\cM$ is a plane which has $\be_z$ tangent to it, and hence $\cM$ is trivially a translator.
\end{proof}

\section{The drift heat equation}
It is well known that the functions $1, \theta$ and the coordinate
functions $x_i$ all
satisfy the heat equation along the mean curvature flow. Along the
rescaled flow we  instead consider rescaled coordinate functions as follows.
\begin{definition}\label{defn:rescaledcoord}
 For any coordinate function $x_i$ on $\CC^n$ we have the \emph{rescaled coordinate function}
$$\tilde{x}_i=e^{-\tau/2}x_i$$
along the rescaled flow $M_{\tau}$.  In particular, we have the \emph{rescaled height}
$\tilde w = e^{-\tau/2} w.$
\end{definition}

Using the above definition, the next result, which is key for our purposes,   follows from a straightforward rescaling.
\begin{lemma}\label{lem:evol.coords} The functions $1, \theta$ and the
  rescaled coordinate functions $\tilde{x}_i$ satisfy the  drift heat
  equation 
  \begin{equation}\label{eq:driftheat} \frac{\partial f}{\partial\tau} = \mathcal{L}_0 f\end{equation}
  along the rescaled flow, where
  \begin{equation}\label{eq:driftL} \mathcal{L}_0 f := \Delta f - \frac{1}{2}\langle \bx, \nabla
    f\rangle \end{equation}
  is the drift Laplacian. 
\end{lemma}

\noindent Note that, when
computing derivatives $\frac{\partial f}{\partial\tau}$, the
rescaled flow has velocity $\mathbf{H} + \frac{1}{2}\mathbf{x}^\perp$.

We will compare solutions of the drift heat equation along the
rescaled flow with solutions on the blow-downs. By
Proposition~\ref{prop:blowdown.line} all possible blow-downs are 
unions $P_1'\cup P_2'$ of two $n$-dimensional subspaces of $\CC^n$ meeting along a
subspace of dimension less than $n$. We therefore study solutions of
\eqref{eq:driftheat} on Euclidean spaces. 

On an $n$-dimensional space $P = \mathbb{R}^n$ we define the drift heat
equation and  drift Laplacian by  \eqref{eq:driftheat} and
\eqref{eq:driftL}. For a solution $f(\bx, \tau)$ of the drift heat equation
on $P$, we define the weighted norm $\Vert f\Vert_{\tau}$ by
\begin{equation}\label{eq:Idefn}
  \Vert f\Vert_\tau^2 = \int_P f(\bx, \tau)^2 e^{-|\mathbf{x}|^2/4}.
  \end{equation}
By \cite[Theorem 0.6]{ColdingMinicozzi:frequency} the function $\log
\Vert f\Vert_\tau$ is convex in $\tau$, and it is linear if and only if $f$
is homogeneous, i.e.~$f(\bx,
\tau) = e^{-\lambda \tau} h(\bx)$, where $h$ is an eigenfunction of
$\mathcal{L}_0$ with eigenvalue $\lambda$. In this case $\log \Vert
f\Vert_\tau^2 = -2\lambda\tau + \log \Vert h\Vert^2$  and we say that $f$
has degree $2\lambda$. The eigenfunctions of the
Ornstein--Uhlenbeck operator $\mathcal{L} = \Delta - \mathbf{x}\cdot\nabla$ on
Euclidean space are well-studied, see e.g.~Bogachev~\cite[Chapter
1]{Bogachev}. The eigenvalues of $\mathcal{L}$ are non-negative integers $k$, and the
corresponding eigenfunctions are degree $k$ homogeneous polynomials
given by products of Hermite polynomials. If $H_k$ is an eigenfunction
of $\mathcal{L}$ with eigenvalue $k$, then the
function $h_k(\bx) = H_k(\bx/\sqrt{2})$ is an eigenfunction of
$\mathcal{L}_0$ with eigenvalue $k/2$. This leads to the following.

\begin{lemma}
\label{lem:plane.drift}
Let $P = \mathbb{R}^n$ and let $x_i$ be coordinate functions on
$P$. The eigenvalues of $\mathcal{L}_0$ on $P$ are given by
non-negative half integers, and so the homogeneous solutions of the
drift heat equation on $P$ have non-negative integer degrees. The
homogeneous solutions with degree 0 are the constants, while those
with degree 1 are spanned by the rescaled coordinate functions $e^{-\tau/2} x_i$.
\end{lemma}

We will be interested in solutions to the drift heat equation on the blow-downs $P_1'\cup P_2'$,
where $P_j'$ are two distinct $n$-dimensional subspaces of $\CC^n$. We
define these as follows.
\begin{definition}\label{defn:pairofsoln} A solution of the (drift) heat equation on
  $P_1'\cup P_2'$ is a pair $u = (u_1, u_2)$, where $u_j$ is a
  solution of the (drift) heat equation on $P_j'$. We define the
  weighted norm $\Vert u\Vert_{\tau}$ of $u$ by $\Vert u\Vert_\tau^2 = \Vert u_1\Vert_\tau^2
  + \Vert u_2\Vert_\tau^2$.  
\end{definition}

We observe that the function $\theta$, equal to the constant
$\overline{\theta}_j$ on $P_j'$, is
a solution of the (drift) heat equation on $P_1'\cup P_2'$ in this sense. Note that we can see $u = (u_1, u_2)$ as one solution to the (drift) heat equation on the (immersed) shrinker $P_1'\cup P_2'$, so we still have by  \cite[Theorem 0.6]{ColdingMinicozzi:frequency}  that
$\log \|u\|_\tau$ is convex, and it is linear if and only if $u$ is homogeneous.

 Lemma~\ref{lem:plane.drift} implies the following.
\begin{lemma}\label{lem:blowdownhom}
  On any blow-down $P_1' \cup P_2'$ the homogeneous solutions of the
  drift heat equation have non-negative integer degrees.
\end{lemma}

Recall our basic assumption that one blow-down is given by
$P_1\cup P_2$, where $P_1\cap P_2=\ell$ is a line. Recall the coordinates $x_1,\ldots,x_{2n-2},z,w$ as chosen in $\S$\ref{ss:setup}, where the
coordinate along $\ell$ is $z$  and $w$ is the height function
vanishing along $P_1\cup P_2$.  
We then have the
following, which also uses
the assumption that the Lagrangian angles of $P_1$ and $P_2$ are
different.
\begin{lemma}\label{lem:P1P2soln}
  The degree 0 solutions of the drift heat equation on $P_1\cup P_2$
  are spanned by $1, \theta$. The degree 1 solutions are spanned by
  $e^{-\tau/2} x_1,\ldots, e^{-\tau/2}x_{2n-2}$ and $e^{-\tau/2} z,
  e^{-\tau/2}z \theta $. 
\end{lemma}
\begin{proof}
  The degree 0 solutions on $P_1\cup P_2$
  are given by  pairs $(c_1, c_2)$ of constants. These are spanned by
  the functions $1, \theta$ since $\theta$ equals two distinct
  constants $\overline{\theta}_j$ on the subspaces $P_j$. 

  The degree 1 solutions on $P_1\cup P_2$ are given by pairs $(f_1,
  f_2)$ of linear functions on $\CC^n$ restricted to the subspaces. According to our choice of
  coordinates in $\S$\ref{ss:setup}, $f_1$ is in the span of $x_1,\ldots, x_{n-1}, z$
  and $f_2$ is in the span of $x_n, \ldots, x_{2n-2}, z$. Since
  $x_1,\ldots, x_{n-1}$ vanish on $P_2$, and $x_n,\ldots, x_{2n-2}$
  vanish on $P_1$, the collection of functions $x_1,\ldots, x_{2n-2}$
  on $P_1\cup P_2$ define the pairs $(x_i,0)$ and $(0,x_j)$, where
  $1\leq i\leq n-1$ and $n\leq j\leq 2n-2$. At the same time $z,
  z\theta $ contain the pairs $(z,0)$ and $(0,z)$ in their span (again since $\theta$ takes distinct values on $P_1,P_2$). 
\end{proof}

\subsection{Limits of solutions of the heat
  equation}\label{ss:heat.limits}  In this subsection we show that if
we are given a solution $u$ of the heat equation along the ancient
mean curvature flow $\cM$, and a sequence of rescalings of $\cM$
converging to a blow-down $P_1' \cup P_2'$ given by a union of distinct
$n$-dimensional subspaces,
then along a subsequence we can extract a normalized limit of $u$,
determining a solution of the heat equation on both planes
separately.

Standard methods allow us to extract limits on compact
sets away from the intersection $E = P_1'\cap P_2'$, and in the limit
we obtain solutions of the heat equation on $P_j' \setminus E$ for
$j=1,2$, which are in $L^\infty$ across $E$.
The main difficulty is that
$E$ may have codimension 1 in $P_j'$ and codimension 1
sets are not removable for solutions of the heat equation. (Consider
for instance the solution given by two different constants on the
lower and upper half planes.) To overcome
this issue it is crucial that the angle $\theta$ differs on the
two subspaces (which we have by Proposition~\ref{prop:blowdown.line}), while at the same time the space-time integral of $|\nabla\theta|^2$
converges to zero as we approach the blow-down. This allows us to show
that the solutions that we obtain in the limit on $P_j' \setminus E$ are
distributional solutions across $E$, and hence smooth. 

To state the result, let $L_t^i$ be a sequence of smooth
solutions of LMCF in $\CC^n$ defined for $t\in [-1,0]$. We
assume that the $L_t^i$ have Euclidean area growth  and uniformly bounded
Lagrangian angles. We assume that $L_t^i\tweak P_1'\cup P_2'$ weakly as
$i\to\infty$, where as above $P_j'$ are $n$-dimensional subspaces meeting along
a subspace $E$ of dimension at most $n-1$. Here, as usual, we view $P_1'\cup
P_2'$ as a static  flow. For the definition of functions with polynomial growth we refer the reader to Definition \ref{def:polynomial_growth}. 

\begin{proposition}
  \label{prop:limit.heat} In the setting above, for each $i$,
  let $u_i$ be a solution to the heat equation on $L_t^i$ for
  $t\in[-1,0]$, with at most polynomial growth.
  Assume further that there is a uniform $C>0$ so that
\begin{equation}\label{eq:limit.heat.bound} 
\int_{L^i_{-1}} u_i^2 e^{-|\bx|^2/4} < C. 
\end{equation}
Then, after passing to a subsequence, we have 
$u_i\to \overline{u}$ where $\overline{u} = (\overline{u}_1,
\overline{u}_2)$ is a solution of the
heat equation on the union $P_1'\cup P_2'$ for $t\in (-1,0]$ in the sense of
Definition~\ref{defn:pairofsoln}. The convergence $u_i \to
\overline{u}$ here means smooth convergence on compact subsets of
$(-1,0] \times \CC^n\setminus E$, i.e.~on compact subsets away from
$t=-1$ and away from the intersection $P_1'\cap P_2'$.
\end{proposition}

\begin{proof}
  We have
  \[ (\partial_t - \Delta) u_i^2 = -2|\nabla u_i|^2. \]
 We apply the monotonicity formula \eqref{eq:Huisken}  (see Proposition \ref{prop:monotonicity}) to $u_i^2$ with different
 centers $(\bx_0, t_0)$ in $\CC^n\times (-1,0]$. 
 Using the uniform bound \eqref{eq:limit.heat.bound} we find 
 that for any $R > 0$ there is a constant $C_R>0$  so that
\begin{equation}\label{eq:lim.bounds}
 \begin{aligned} \sup_{B_R(0)\times [-1+R^{-1},0]} |u_i| &< C_R, \\
      \int_{-1+R^{-1}}^0 \int_{B_R(0)\cap L^i_t} |\nabla u_i|^2 &<
      C_R. \end{aligned} 
      \end{equation}
      
Let $\theta_i$ be the Lagrangian angle on $L^i_{-1}$ and
$\overline{\theta}_1\not= \overline{\theta}_2$ the (constant) Lagrangian angles
    on $P_1',P_2'$. 
As in \cite[Theorem A]{Neves.ZM}, we have that for 
    all $s\in (-1,0)$, $f\in C^2(\RR)$ and compactly supported smooth functions $\phi$, 
    \begin{equation*}
      \lim_{i\to\infty} \int_{L^i_s} f(\theta_i) \phi = \sum_{j=1}^2 \int_{P_j}
      f(\overline{\theta}_j)\phi.
    \end{equation*}
    Since $\overline{\theta}_1\not=\overline{\theta}_2$, we
    can choose $f\in C^2(\RR)$ such that $f(\overline{\theta}_1)=1$ and
    $f(\overline{\theta}_2)=0$, and we fix such a function $f$ for the rest of the proof.
    
       We also  fix a smooth function $\chi$ compactly supported in $B_R(0)\times
    (-1,0)$.  Then we have
    \begin{equation}\label{eq:lim.evol} 
    \begin{aligned}
        \frac{d}{dt} \int_{L^i_t} f(\theta_i) u_i \chi &= \int_{L^i_t}
        f'(\theta_i) (\Delta \theta_i) u_i \chi + \int_{L^i_t} f(\theta_i)
        (\Delta u_i) \chi + \int_{L^i_t} f(\theta_i) u_i \partial_t\chi \\
        &\quad - \int_{L^i_t} f(\theta_i) u_i \langle J\nabla\theta_i, D\chi\rangle -
        \int_{L^i_t} f(\theta_i) u_i \chi |\nabla\theta_i|^2,\\
      \end{aligned} 
      \end{equation}
      where $D$ denotes the ambient derivative on  Euclidean space,
      using the fact that both $u_i$ and $\theta_i$ solve the heat
      equation on $L^i_t$ and $\bH=J\nabla\theta_i$. Note that, since
      $\chi$ has compact support, we may use \eqref{eq:lim.bounds} and
      the fact that the (spacetime) $L^2$-norm of $|\bH| =|\nabla \theta_i|$ goes to
    zero as $i\to\infty$ (see \cite[Lemma 5.4]{Neves.ZM})  to deduce that
\begin{equation*}
\begin{split}
 \int_{-1}^0\int_{L^i_t}
        f'(\theta_i) (\Delta \theta_i) u_i \chi\rd t &- \int_{-1}^0\int_{L^i_t} f(\theta_i) u_i \langle J\nabla\theta_i, D\chi\rangle\rd t \\
        &-
        \int_{-1}^0\int_{L^i_t} f(\theta_i) u_i \chi |\nabla\theta_i|^2\rd t\to 0\qquad\text{as $i\to\infty$.}
\end{split}
\end{equation*}

Therefore, since $\chi$ has compact support in $B_R(0)\times (-1,0)$, if we integrate \eqref{eq:lim.evol} with respect to $t$ on $[-1,0]$, we have that
    \begin{equation}\label{eq:lim.epsi.1} \int_{-1}^0 \int_{L^i_t} f(\theta_i) u_i \partial_t \chi\rd t =
      -\int_{-1}^0 \int_{L^i_t} f(\theta_i) (\Delta u_i) \chi\rd t  +
      \epsilon_i,
      \end{equation}
    where $\epsilon_i\to 0$ as $i\to\infty$.  (Note that $\epsilon_i$ will
    depend on $\chi$.)
    Integrating by parts on the right-hand side of \eqref{eq:lim.epsi.1} we get
\begin{equation}\label{eq:lim.epsi.2}
\begin{split}
\int_{-1}^0 \int_{L^i_t} f(\theta_i) u_i \partial_t \chi\rd t &= 
      \int_{-1}^0 \int_{L^i_t} f(\theta_i) \langle\nabla u_i,\nabla \chi\rangle\rd t\\
    &\quad  +\int_{-1}^0 \int_{L^i_t} f'(\theta_i) \chi \langle\nabla u_i, \nabla \theta_i\rangle\rd t
      + \epsilon_i. 
      \end{split}
      \end{equation}
 Again from \eqref{eq:lim.bounds} and the fact $\nabla\theta_i$ converges
    to zero in $L^2$ (in spacetime), we may absorb the second integral on the right-hand side of \eqref{eq:lim.epsi.2} into $\epsilon_i$.
    We can then integrate by parts in the first term on the right-hand side of \eqref{eq:lim.epsi.2} and absorb
    another term involving $\nabla \theta_i$ by the same argument into $\epsilon_i$ to get
    \begin{equation}\label{eq:lim.epsi.3} \int_{-1}^0 \int_{L^i_t} f(\theta_i) u_i \partial_t \chi\rd t = 
      -\int_{-1}^0 \int_{L^i_t} f(\theta_i) u_i \Delta \chi\rd t +
      \epsilon_i.
      \end{equation}
      Recall that $f(\overline{\theta}_1)=1$ and $f(\overline{\theta}_2)=0$.  
  Since the $u_i$ are uniformly bounded on the support of $\chi$ by \eqref{eq:lim.bounds}, and we have good
    convergence away from the line $\ell=P_1\cap P_2$, the contribution as we pass to the limit as $i\to\infty$ in \eqref{eq:lim.epsi.3} near the singular set $\ell$ is negligible. Therefore, we can
    pass to the limit in \eqref{eq:lim.epsi.3} along a subsequence, and get that the subsequential limit $\overline{u}_1$ of the $u_i$ on $P_1$ satisfies
    \[ \int_{-1}^0 \int_{P_1} \overline{u}_1 \partial_t \chi\rd t = -\int_{-1}^0
      \int_{P_1} \overline{u}_1\Delta \chi\rd t. \]
     This means that the limit $\overline{u}_1$ is a bounded distributional solution of the   heat equation on $P_1$ so it follows that $\overline{u}_1$ is a classical solution on $P_1$.  
     
 Repeating the argument starting with the subsequence converging to $\overline{u}_1$ on $P_1$ and changing the choice of function $f$ so that it takes the value $1$ on $\overline{\theta}_2$ and $0$ on $\overline{\theta}_1$ yields the result.
\end{proof}

Since the
  drift heat equation and usual heat equation are related by
  rescaling, one can apply Proposition \ref{prop:limit.heat} to
  sequences of solutions of the drift heat equation along rescaled
  mean curvature flows. In
  particular, suppose that we have a sequence of rescaled flows
  $M^i_{\tau}$, for $\tau\in[-1,0]$, converging to $P_1'\cup P_2'$
  weakly. Recall  the weighted $L^2$-norm defined in \eqref{eq:Idefn} and let $u_i$ be solutions of the drift heat equation on
  $M^i_\tau$, with $\Vert u_i\Vert_{-1} \leq
  1$, and such that the $u_i$ have polynomial growth.   
  Proposition~\ref{prop:limit.heat} implies that after passing to
  a subsequence we have $u_i\to \overline{u}$, for a solution
  $\overline{u}$ of the drift heat equation on $P_1'\cup P_2'$ for
  $\tau\in (-1,0]$. We have the following additional information,
  saying that for $\tau > -1$ the weighted $L^2$-norms of the $u_i$ cannot concentrate
  near $E$ and near infinity. 
  \begin{lemma}\label{lem:L2conv}
    Under the setup above we have
    \begin{equation}\label{eq:L2conv}
      \Vert \overline{u}\Vert_{\tau} = \lim_{i\to\infty} \Vert
      u_i\Vert_{\tau} \leq \liminf_{i\to\infty} \Vert u_i\Vert_{-1}, \end{equation}
    for $\tau\in (-1,0]$. 
  \end{lemma}
  \begin{proof}
    The inequality $\Vert u_i\Vert_\tau \leq \Vert u_i\Vert_{-1}$ for
    $\tau > -1$ 
    follows immediately from the monotonicity
    formula~\eqref{eq:Huisken} and the observation that $(\partial_t - \Delta)u_i^2
    \leq 0$.
    
    For $r, R > 0$ let us write $A_{r,R} = B_R(0) \setminus B_{r}(E)$,
    where $E=P_1'\cap P_2'$ and $B_{r}(E)$ denotes the
    $r$-neighbourhood of $E$. Let $\delta > 0$. 
    From Proposition~\ref{prop:limit.heat} we know that for any $r, R >
    0$, and $\tau\in [-1+\delta, 0]$ we have
    \[ \lim_{i\to\infty} \int_{M^i_\tau \cap A_{r,R}} u_i^2 e^{-|{\bf
        x}|^2/4} = \int_{(P_1'\cup P_2')\cap A_{r,R}} \overline{u}^2\,
      e^{-|{\bf x}|^2/4}. \]
     To prove \eqref{eq:L2conv} it is enough to show that for any
     $\epsilon, \delta > 0$, there are $r, R > 0$ such that for all $i$
     and $\tau\in [-1+\delta, 0]$      we have
     \[ \int_{M^i_\tau \setminus A_{r,R}} u_i^2\, e^{-|{\bf x}|^2/4} <
       \epsilon. \]
     First, using the log Sobolev inequality due to Ecker~\cite[Theorem
     3.4]{Ecker.logSobolev}, we have a $p > 1$, depending on $\delta > 0$, such
     that
     \[ \left(\int_{M^i_\tau} |u_i|^{2p} e^{-|{\bf x}|^2/4} \right)^{1/p}
       < C, \]
     for a uniform $C$, as long as $\tau\in [-1+\delta, 0]$. It
     follows using H\"older's inequality that, given $R > 0$, we have
     \[ \int_{M^i_\tau \setminus B_R(0)} |u_i|^2\, e^{-|{\bf x}|^2/4}
       \leq C\left(\int_{M^i_\tau\setminus B_R(0)} e^{-|{\bf
             x}|^2/4}\right)^{1-1/p}, \]
     and so using the Euclidean area bounds for $M^i_\tau$ we can find
     an $R$ (depending on $\delta, \epsilon$) such that
     \begin{equation}\label{eq:b10}
       \int_{M^i_\tau \setminus B_R(0)} |u_i|^2\, e^{-|{\bf x}|^2/4}
       \leq \frac{\epsilon}{2}, \end{equation}
     for $\tau\in [-1+\delta, 0]$.

     Viewing $R$ (and $\delta$) as fixed, 
     the uniform bound in \eqref{eq:lim.bounds} implies that if $r$ is
     sufficiently small (depending on $\epsilon, \delta, R$), then 
     \[ \int_{M^i_\tau\cap B_R(0)\cap B_{r}(E)} |u_i|^2 e^{-|{\bf
           x}|^2/4} < \frac{\epsilon}{2}. \]
     Combined with \eqref{eq:b10} this implies
     \[ \int_{M^i_\tau\setminus A_{r,R}} |u_i|^2 e^{-|{\bf
           x}|^2/4} < \epsilon, \]
     as required. 
  \end{proof}

\subsection{Three annulus lemma}  \label{sec:three-annulus} A well-known method for controlling the
growth of solutions of PDEs is the three annulus lemma, see for
example \cite{Simon:asymptotics}. In this subsection we prove a version of
the three annulus lemma for solutions of the drift heat equation along
the rescaled flow. We use an argument by contradiction, similar to
Simon~\cite{Simon:asymptotics}, based on the monotonicity of frequency shown by
Colding--Minicozzi~\cite{ColdingMinicozzi:frequency}.
Related ideas are also applied in
 \cite{CMcomplexity}. 

In this subsection we assume that $L_\tau$ is a rescaled Lagrangian
mean curvature flow  such that, along a sequence $\tau_i\to -\infty$, we
have $L_{\tau_i} \tweak P_1\cup P_2$. In addition we assume, as
before, that the $L_\tau$ have uniformly bounded area ratios,   uniformly bounded Lagrangian angle and are almost
calibrated for $n\geq 3$.   


\begin{proposition}
\label{prop:three.annulus.lemma}
 For any $s \not\in \ZZ$ there is a $T_0=T_0(s) > 0$ with the
 following property. Suppose that $u$ is a solution of the drift heat
 equation \eqref{eq:driftheat} on the rescaled flow $M_\tau$ for
 $\tau\in[-T-2,-T]$ with $T > T_0$, such that $u$ has   polynomial
 growth. 
 If in addition we have that the weighted $L^2$-norm defined in
 \eqref{eq:Idefn} satisfies 
 \[  \| u\|_{-T-1}\geq e^{s/2}\|u\|_{-T},\]
 then we also have
 \[ \|u\|_{-T-2}\geq e^{s/2}\|u\|_{-T-1}.\]
\end{proposition}
\begin{proof}
  We argue by contradiction. Suppose that there is a sequence of
  solutions $u_i$ to \eqref{eq:driftheat} on intervals $[-T_i-2,
  -T_i]$ with $T_i\to \infty$, such
  that
  \begin{equation}\label{eq:ui1}
     \|u_i \|_{-T_i-1} \geq e^{s/2} \| u_i\|_{-T_i}, \end{equation}
  but
  \begin{equation} \label{eq:ui2} \|u_i \|_{-T_i-2} < e^{s/2} \|
    u_i\|_{-T_i-1}. \end{equation} 
  By rescaling we can assume that $\| u_i\|_{-T_i-1} = 1$ for all
  $i$. It follows from \eqref{eq:ui2} that then $\| u_i \|_{-T_i-2} < e^{s/2}$. We can
  apply Proposition~\ref{prop:limit.heat} to time translations of the
  $u_i$, and along a subsequence we can extract a limit $\overline{u}$
  satisfying the drift heat equation along a blow-down $P_1'\cup P_2'$
  of the flow $L_\tau$ on the interval $(-2,0]$. Using \eqref{eq:L2conv} we have
  \begin{equation} \label{eq:ubar0} \Vert \overline{u}\Vert_{-1} = 1,
  \end{equation}
  and at the same time the local uniform convergence of $u_i$ to
  $\overline{u}$, together with \eqref{eq:ui1} and \eqref{eq:ui2}, implies
  \begin{equation} \label{eq:ubar1} \Vert \overline{u}\Vert_0 \leq e^{-s/2}, \qquad \Vert
    \overline{u} \Vert_{\tau} \leq e^{s/2} \text{ for all
    }\tau\in(-2,0].
  \end{equation}
  By \cite[Theorem 0.6]{ColdingMinicozzi:frequency} we know that $\log
  \Vert \overline{u}\Vert^2_\tau$ is convex. From \eqref{eq:ubar0} and
  \eqref{eq:ubar1} it follows that 
  $\log \Vert \overline{u}\Vert_\tau^2$ is linear with slope $s$. By  \cite[Theorem
  0.6]{ColdingMinicozzi:frequency}  
  $\overline{u}$ must be homogeneous with degree $s$. 
  By Lemma~\ref{lem:blowdownhom} the
  homogeneous solutions on any blow-down have integer degrees, so since
  $s\not\in\ZZ$, we have a contradiction.  
\end{proof}

We can use the three annulus lemma to extract the leading order
behaviour of ancient solutions to the heat equation as follows.
\begin{proposition}\label{prop:homlimit}
  Suppose that $u$ is a non-zero solution of the drift heat equation along the
  rescaled flow
  $M_\tau$ for $\tau\in (-\infty, 0]$, with  polynomial growth. Suppose
  that for some $C, d > 0$ we have $\Vert u\Vert_{\tau}^2 \leq C
  e^{-d\tau}$ for all $\tau < 0$. Let $\tau_i\to -\infty$ be integers. Up to
  choosing a subsequence we have the following. The translated (rescaled) flows
  $L^i_\tau = M_{\tau-\tau_i}$ converge weakly to a blow-down
  $P_1'\cup P_2'$, and the normalized translated solutions
  \[ u_i(\tau) = \Vert u \Vert_{\tau_i}^{-1}\, u(\tau-\tau_i) \]
  converge to a non-zero homogeneous solution $\overline{u}$ of the
  drift heat equation on $P_1'\cup P_2'$ for $\tau\in [-2,0]$. The
  convergence is locally smooth on $[-2,0]$ away from $P_1'\cap
  P_2'$, and also in $L^2$ as in \eqref{eq:L2conv}. 
\end{proposition}
\begin{proof}
  Let $s_0 > d$ for some $s_0\not\in\ZZ$. We claim that
  there is a $\tau_0 < 0$ such that we then have
  \begin{equation} \label{eq:us0} \Vert u\Vert_{\tau-1} \leq e^{s_0/2}
    \Vert u\Vert_{\tau} \end{equation} 
  for all $\tau < \tau_0$. If this were not the case, then
  Proposition~\ref{prop:three.annulus.lemma} would imply that in fact
  $\Vert u\Vert_{\tau-k} \geq e^{s_0/2} \Vert u\Vert_{\tau-k+1}$ for all
  integers $k > 0$ and some $\tau$, but this would eventually contradict the
  assumption $\Vert u\Vert_{\tau-k}^2 \leq C e^{-d(\tau-k)}$.

  The growth condition \eqref{eq:us0} together with the normalization
  of $u_i$ implies that $\Vert u_i\Vert_{-3} \leq e^{3s_0/2}$. Using
  Proposition~\ref{prop:limit.heat} we can extract a limit
  $\overline{u}$ along a subsequence
  on $P_1'\cup P_2'$. The convergence is locally smooth
  on $(-3,0]$ away from $P_1'\cap P_2'$, and using
  Lemma~\ref{lem:L2conv} the convergence is in $L^2$ for
  $\tau\in[-2,0]$ as required.

  It remains to argue that $\overline{u}$ is homogeneous. For this
  note that Proposition~\ref{prop:three.annulus.lemma} implies that
  for any $s\not\in\ZZ$ one of the following must hold:
  \begin{enumerate}
  \item $\Vert u\Vert_{\tau-1} \geq e^{s/2}\Vert u\Vert_{\tau}$ for
    all sufficiently negative integers $\tau$,
  \item $\Vert u\Vert_{\tau-1} \leq e^{s/2}\Vert u\Vert_{\tau}$ for all
      sufficiently negative integers $\tau$, 
    \end{enumerate}
    since if (1) holds for some sufficiently negative $\tau$ then it
    must hold for $\tau-k$ for all integers $k > 0$ by
    Proposition~\ref{prop:three.annulus.lemma} as well. It follows that there is
    some $s_1\in \RR$ such that (1) holds for all $s < s_1$, and (2)
    holds for all $s > s_1$. We deduce that in the limit we
    have
    \[ \Vert \overline{u}\Vert_{-2} = e^{s_1/2} \Vert
      \overline{u}\Vert_{-1}, \quad \Vert \overline{u}\Vert_{-1} = e^{s_1/2} \Vert
      \overline{u}\Vert_{0}. \]
    The convexity of $\log \Vert\overline{u}\Vert_\tau$ then implies
    that $\log \Vert \overline{u}\Vert_{\tau}$ is linear, from which it follows that $\overline{u}$ is homogeneous. 
\end{proof}

We will also need the following variant of the three annulus lemma,
similar to Donaldson--Sun~\cite[Proposition 3.11]{DS2}. 
\begin{proposition}\label{prop:3ann2}
  Let  $V_{\leq 1}$ be the space of solutions of the drift
  heat equation along $M_\tau$ given by the span of $1, \theta$ and
  $e^{-\tau/2}x_i$ for the coordinate functions $x_i$. Let $V\subset V_{\leq
    1}$ be any subspace  and let $u$ be a solution of the drift heat
  equation along $M_\tau$ with polynomial growth. Suppose that
  there is a constant $C > 0$ such that
  \begin{equation}\label{eq:u32}
    \Vert u\Vert_{\tau}^2 \leq C e^{-3\tau/2}, \text{ for all } \tau <
    -1. \end{equation}

  For any $\tau$ let
  $\Pi_\tau u := u - f$, where $f\in V$ and $u-f$ is
  orthogonal to $V$ at time $\tau$:
  \[ \langle u - f, g\rangle_\tau := \int_{M_\tau} (u-f)g\, e^{-|{\bf
        x}|^2/4} = 0, \text{ for all } g\in V. \]
  Given $s\not\in
  \ZZ$, there is a $T_0 > 0$ with the following property. If
  \[ \Vert \Pi_{-T-1}u \Vert_{-T-1} \geq e^{s/2} \Vert \Pi_{-T}u
    \Vert_{-T} \]
  for some $T > T_0$, then
  \[ \Vert \Pi_{-T-2} u\Vert_{-T-2} \geq e^{s/2} \Vert
    \Pi_{-T-1}u\Vert_{-T-1}. \]
\end{proposition}
\begin{proof}
  The proof is by contradiction, similar to that of
  Proposition~\ref{prop:three.annulus.lemma}. Suppose that we have a
  sequence $T_i\to\infty$  and corresponding $u_i$ such that
  \begin{equation}\label{eq:ui11} \Vert \Pi_{-T_i-1}u_i \Vert_{-T_i-1} \geq e^{s/2} \Vert \Pi_{-T_i}u_i
    \Vert_{-T_i}, \end{equation}
  and at the same time
  \begin{equation} \label{eq:ui12} \Vert \Pi_{-T_i-2} u_i\Vert_{-T_i-2} < e^{s/2} \Vert
    \Pi_{-T_i-1}u_i\Vert_{-T_i-1}. \end{equation}
   Let  $v_i=\Pi_{-T_i-2}u_i$, so that $v_i$ is orthogonal to $V$ at
   $\tau=-T_i-2$. By scaling we can assume that $\Vert v_i\Vert_{-T_i-1}
   =1$. It follows that $\Vert \Pi_{-T_i-1} v_i\Vert_{-T_i-1}\leq 1$,
   and so by \eqref{eq:ui12} we have $\Vert v_i\Vert_{-T_i-2}\leq
   e^{s/2}$.  We claim that for sufficiently large $i$ we also
   have
   \begin{equation}\label{eq:vigrowth}
     \Vert v_i\Vert_{-T_i-3} \leq e^{4/5} \Vert v_i\Vert_{-T_i-2}.
   \end{equation}
   If \eqref{eq:vigrowth} did not hold, Proposition~\ref{prop:three.annulus.lemma} would
   imply that for some constant $C > 0$ and for
   all integers $k > 3$ we would have $\Vert
   v_i\Vert_{-T_i-k} \geq C^{-1} e^{4k/5}$. At the same time $v_i = u_i -
   f_i$ for some $f_i\in V$, and since both $u_i$ and $f_i$ satisfy an
   estimate of the form \eqref{eq:u32}, we get a contradiction. Thus
   \eqref{eq:vigrowth} holds, and so we have a uniform bound $\Vert
   v_i\Vert_{-T_i-3} \leq e^{s/2 + 4/5}$.

 Applying Proposition~\ref{prop:limit.heat} we have that, along a
   subsequence and after time translations, the $v_i$ converge to a limit solution
   $\overline{v}$ of the drift heat equation on a blow-down $P_1'\cup
   P_2'$ for $\tau\in (-3,0]$. It follows using \eqref{eq:L2conv} that 
   $\Vert \overline{v}\Vert_{-2}\leq e^{s/2}$ and $\Vert
   \overline{v}\Vert_{-1}=1$.   

   We claim that we also have 
\begin{equation}\label{eq:ubar20} 
\Vert \overline{v}\Vert_0 \leq
   e^{-s/2},\end{equation}
 in which case we will reach a contradiction just like in
   the proof of Proposition~\ref{prop:three.annulus.lemma}. Note that
   the new difficulty is that we only have the bound $\Vert \Pi_{-T_i}
   v_i\Vert_{-T_i}\leq e^{-s/2}$, and the norm of $v_i$ can be larger than
   that of its projection $\Pi_{-T_i}v_i$. 

   To see that \eqref{eq:ubar20} holds we show that under our
   assumption that $v_i$ is orthogonal to $V$ at time $-T_i-2$, we
   have that $v_i$ is also approximately orthogonal to $V$ at time
   $-T_i$. Let $g\in V$ and consider normalizations $g_i$ of $g$ such
   that $\Vert g_i\Vert_{-T_i} = 1$. By
   Proposition~\ref{prop:homlimit}, after taking a further subsequence
   we can assume that the $g_i$ converge to a homogeneous limit
   $\overline{g}$ on $P_1'\cup P_2'$, on the time
   interval $[-2, 0]$, satisfying the drift heat equation. We can
   apply the $L^2$-convergence \eqref{eq:L2conv} to $v_i \pm g_i$, together with our
   assumption $\langle v_i, g_i\rangle_{-T_i-2}=0$ to find that $\langle
   \overline{v}, \overline{g}\rangle_{-2} = 0$. Since $\overline{g}$
   is homogeneous, this implies that $\langle \overline{v},
   \overline{g}\rangle_\tau = 0$ for all $\tau\in [-2,0]$. It follows
   from the $L^2$-convergence we have $\langle v_i, g_i\rangle_{-T_i}
   \to 0$. Since this applies to all $g\in V$, we find that
   \[ \lim_{i\to \infty} \frac{ \Vert \Pi_{-T_i} v_i\Vert_{-T_i}}{
       \Vert v_i\Vert_{-T_i}} = 1, \]
   and it follows that $\Vert \overline{v}\Vert_{0} \leq
   e^{-s/2}$. This leads to a contradiction as discussed above. 
 \end{proof}

Let us use coordinates $x_1,\ldots, x_{2n-2}, z, w$ as in $\S$\ref{ss:setup}.
Recall that a blow-down of our ancient
rescaled flow $M_\tau$ along a sequence of scales $\tau_i\to -\infty$
is given by $P_1\cup P_2$, where $P_1\cap P_2=\ell$ is a
line, the coordinate $w$ vanishes on $P_1\cup P_2$ and $J\nabla z=
\nabla w$. Without loss of generality we can assume that the $\tau_i$
are all integers. Let us
write $L^i_t$ for the corresponding sequence of flows for $t\in
[-2,0)$ converging weakly to $P_1\cup P_2$. We then have
the following dichotomy.

\begin{proposition}\label{prop:2cases}
  Either we have that $w = a + b\theta$ for some constants $a,b$ along
  our flow $L_t$, so $L_t$ is a translator, or up to choosing a subsequence
  of the $\tau_i$ we can find a sequence of linear functions
  $\phi_i \in \Span\{x_1,\ldots,x_{2n-2}, z\}$ with $\phi_i\to 0$ and a sequence
  $\sigma_i\to 0$ such that along the sequence $L^i_t$ converging to
  $P_1\cup P_2$ we have
  \[ \sigma_i^{-1} (w - \phi_i) \to z\theta \quad\text{as $i\to\infty$,} \]
  where the convergence is in $L^2$ and locally uniformly away from
  the line $\ell$. 
\end{proposition}
\begin{proof} Recall Definition \ref{defn:rescaledcoord} and let $V=\Span\{1, \theta, \tilde{x}_1,\ldots,\tilde{x}_{2n-2}, \tilde{z}\}$. 

Suppose first that the rescaled height $\tilde{w}$ is in  $V$. Note that $M_{\tau_i} \tweak
P_1\cup P_2$ and $w$ vanishes on $P_1\cup P_2$, but non-trivial linear
combinations of $x_1,\ldots,x_{2n-2}, z$ do not
vanish on $P_1\cup P_2$. This implies that we must  have
$\tilde{w} = a + b\theta$ for constants $a, b$. By
Proposition~\ref{prop:translator.height} the flow $L_t$ is a
translator. 

Suppose now that $\tilde{w}$ is not in $V$. 
We apply Proposition~\ref{prop:3ann2} to $\tilde{w}$
along the flow  with $V$
as chosen. For any integer $k <
0$ let us write
\[ \tilde{w}_k = \frac{ \Pi_k \tilde{w} }{\Vert \Pi_k
    \tilde{w}\Vert_k}. \]
Note that by our assumption $\Pi_k \tilde{w} \not=0$ for all $k$. 
Using Proposition~\ref{prop:3ann2} together with the argument in the
proof of Proposition~\ref{prop:homlimit} we find that along a
subsequence $k_i = \tau_i \to -\infty$, time translations of the
$\tilde{w}_{k_i}$ converge to a homogeneous solution $\overline{w}$
of the drift heat equation on
$P_1\cup P_2$, which is orthogonal to the solutions $1, \theta,
\tilde{x}_1,\ldots,\tilde{x}_{2n-2}, \tilde{z}$. At the same time the growth rate of
$\overline{w}$ can be at most degree 1, so by Lemma~\ref{lem:P1P2soln}
we must have $\overline{w} = c e^{-\tau/2}z\theta $ for a non-zero
constant $c$.

To finish the proof we need to consider how the $\tilde{w}_k$ are
related for different $k$. By definition we have
\begin{equation}\label{eq:twkeqn}
  \tilde{w}_k = \gamma_k \tilde{w}_{k+1} + a_k + b_k\theta + F_k,
\end{equation}
where $\gamma_k, a_k, b_k$ are constants, and $F_k\in
\Span\{\tilde{x}_1,\ldots,\tilde{x}_{2n-2}, \tilde{z}\}$. Using
Proposition~\ref{prop:3ann2}, 
and arguing as in the proof of Proposition \ref{prop:homlimit}, we know that for any subsequence $k_j\to
-\infty$ there is a further subsequence along which the
$\tilde{w}_{k_j}$ (translated in time) converge to a homogeneous
solution along some blow-down $P_1'\cup P_2'$, with degree 1 
 which is orthogonal to $V$. Since
$\Vert \tilde{w}_k\Vert_k = 1$ for all $k$, it follows that $\gamma_k
\to e^{-1/2}$ and $\Vert a_k\Vert_k, \Vert b_k\theta\Vert_k, \Vert
F_k\Vert_k \to 0$ as $k\to -\infty$. Note that the norms $\Vert
1\Vert_k, \Vert \theta\Vert_k$ are uniformly bounded away from 0 and
$\infty$ for all $k$, using the fact that on all blow-downs $P_1'\cup P_2'$
the angle $\theta$ equals the same constants $\overline{\theta}_1,
\overline{\theta}_2$ on the two subspaces $P_1', P_2'$. Therefore
$a_k, b_k \to 0$.

Let us
define the constants $\mu_k$ by $\mu_0=1$ and $\gamma_k =
\mu_{k+1}/\mu_k$ for all sufficiently negative integers $k $. 
From \eqref{eq:twkeqn} we have
\[ \mu_k\tilde{w}_k = \mu_{k+1}\tilde{w}_{k+1} + \mu_k(a_k +b_k\theta) + \mu_k F_k, \]
and so
\begin{equation}\label{eq:tweq2} \tilde{w}_k = \mu_k^{-1} \tilde{w}_0 + \mu_k^{-1}\sum_{i=-1}^k \mu_i(a_i + b_i\theta) +
  \mu_k^{-1}\sum_{i=-1}^k \mu_i F_i. 
\end{equation}
Using that $\mu_{k+1}/\mu_k \to e^{-1/2}$ and $a_k, b_k\to 0$, it
follows that 
\begin{equation}\label{eq:tweq3} \left\Vert \mu_k^{-1} \sum_{i=-1}^k \mu_i(a_i +
    b_i\theta)\right\Vert_k \to 0. 
\end{equation}
At the same time, since $\tilde{w}_0$ is the normalized
$L^2$-projection of $w$ orthogonal to $V$ (at time $\tau=0$), we have $\tilde{w}_0 = c_0 \tilde{w} + c_1 + c_2\theta + F$, where $c_0, c_1,
c_2$ are constants with $c_0\not=0$ and $F\in\Span\{
\tilde{x}_1,\ldots,\tilde{x}_{2n-2}, \tilde{z}\}$. Using \eqref{eq:tweq2} and \eqref{eq:tweq3}
we can  write
\[ \mu_k^{-1}c_0( \tilde{w} - \tilde{\phi}_k) = \tilde{w}_k + E_k, \]
where $\Vert E_k\Vert_k \to 0$ and $\tilde{\phi}_k$ is in
the span of $\tilde{x}_1,\ldots,\tilde{x}_{2n-2}, \tilde{z}$. Along our subsequence $k_i$ we
have $\tilde{w}_{k_i} \to \overline{w} = ce^{-\tau/2} z\theta $, and
so as required we obtain a sequence $L^i_t$ converging to $P_1\cup
P_2$,  and $\sigma_i \neq 0$, $\phi_i\in  \Span\{x_1,\ldots,x_{2n-2}, z\}$ satisfying
\[ \sigma_i^{-1}(w - \phi_i) \to z\theta. \]
It remains to show that $\sigma_i, \phi_i \to 0$. Note that since $w$
vanishes on $P_1\cup P_2$, on $L^i_{-1}$ we have $\Vert
w\Vert_{L^i_{-1}} \to 0$ as $i\to\infty$, while $\Vert x_j\Vert_{L^i_{-1}}$ and $\Vert
z\Vert_{L^i_{-1}}$ are bounded away from $0$ and $\infty$. It
follows that if $\phi_i\not\to 0$, along a
subsequence, then also $\sigma_i \not \to 0$ along this subsequence, and we would  have
$\sigma_i^{-1}\phi_i \to z\theta$ in $L^2$, but this contradicts the fact that
on $P_1\cup P_2$ the function $z\theta$ is $L^2$-orthogonal to   $x_1,\ldots,x_{2n-2}, z$. Therefore we must have $\phi_i \to 0$, which
implies that $\Vert w - \phi_i\Vert_{L^i_{-1}}\to 0$ and so $\sigma_i\to 0$ as
well. 
\end{proof}

In the next section we will show using a topological argument
 that the second alternative in Proposition~\ref{prop:2cases} leads to
 a contradiction. This will complete the proof of our main result.

  \section{Linking argument}\label{sec:linking}

 In this section we use a topological argument to rule out the second
 alternative in Proposition~\ref{prop:2cases}. 
Throughout this section we let
$(-\infty,0) \ni t \mapsto L_t\subset \CC^n$ be a smooth, exact, ancient solution of
LMCF with uniformly bounded
area ratios and Lagrangian angle and which is almost calibrated for $n\geq 3$. Recall that for a  positive
sequence $\lambda_i\to 0$, we consider the sequence of parabolically rescaled  flows 
$$(-\infty, 0) \ni t\mapsto L_t^i=\lambda_iL_{\lambda_i^{-2}t}\, . $$
We assume that as $i\to\infty$ the flows $t\mapsto L^i_t$
converge weakly to the static flow $(-\infty, 0) \ni t\mapsto P_1\cup P_2$, where $P_1,
P_2$ are $n$-dimensional Lagrangian subspaces meeting along a line $\ell$. We
write $\overline{\theta}_j$ for the Lagrangian angles of $P_j$ as
before, where $\overline{\theta}_1 = - \overline{\theta}_2$.

Since the $L^i_t$ are exact, they admit primitives $\beta_i$ of the
Liouville form as in Definition~\ref{dfn:exact}. We have the following
(see Neves~\cite[Proposition 6.1]{Neves.ZM}).
\begin{lemma}\label{lem:betaexist}
  We can choose the primitives $\beta_i$ along the flows $L^i_t$ such
  that
  \[ (\partial_t - \Delta) (\beta_i + 2t\theta_i) = 0. \]
\end{lemma}

Since $|\nabla\beta_i| = |\bx^\perp|$ and $L^i_{-1}$
converges to the union $P_1\cup P_2$ locally smoothly away from
$\ell$, we have that $\beta_i|_{L^i_{-1}} \to\overline{\beta}_j$ as
$i\to\infty$ locally smoothly on each plane $P_j$ away from $\ell$, for suitable
constants $\overline{\beta}_j$. Similarly
$\theta_i\to\overline{\theta}_j$ locally smoothly on $P_j$ away from
$\ell$, as $i\to\infty$.  Given this, we make the following
definition. 

\begin{definition}\label{dfn:scale.cpts}
Since the $L_t^i$ are exact, and almost calibrated for $n\geq 3$, by \cite[Theorem
4.2]{Neves:survey}  there exists a set $\mathcal{E}\subseteq
(-2,0)$ of measure zero so that whenever
$s'\in(-2,0)\setminus\mathcal{E}$, we have two distinct connected
components $\Sigma_{1,s'}^i,\Sigma_{2,s'}^i$ of $B_3(0)\cap L^i_{s'}$
(after possibly passing to a subsequence) intersecting $B_2(0)$ and
converging  (as Radon measures) to the planes $P_1,P_2$ respectively
in $B_2(0)$. 
 Note that there might be more connected components of $B_3(0)\cap L^i_{s'}$, but the components $\Sigma_{1,s'}^i,\Sigma_{2,s'}^i$ are uniquely determined, and the remaining components converge to zero as Radon measures. 

Let $s_1\in(-1/2,0)\setminus\mathcal{E}$  so that
$\overline{b}_1\not= \overline{b}_2$, where
\begin{equation}\label{eq:distinct}
\overline{b}_j:=
\cos(\overline{\beta}_j-2(1+s_1)\overline{\theta}_j). 
\end{equation}
This is always possible since
$\overline{\theta}_1\neq\overline{\theta}_2$ and $\mathcal{E}$ has
measure $0$. We then let $\Sigma^i_j=\Sigma^i_{j,s_1}\cap B_2(0)$ in the notation
above.   
\end{definition}
 
\subsection{Approximate solutions of the heat equation} We now
prove our first key result, which provides solutions of the heat
equation, which on the two components $\Sigma^i_j$, for $j=1,2$, approximate
$\overline{b}_jz$ pointwise. 
 
 \begin{proposition}\label{prop:Biz.approx.heat}  Recall the notation
   of Definition \ref{dfn:scale.cpts}. Let  
\begin{equation}\label{eq:Bi}
B_i= \cos(\beta_i+2(t-s_1)\theta_i)
\end{equation}
and let $h_i$ be the solution of the heat equation along $L^i_t$ with polynomial growth (see Proposition \ref{prop:caloric_unique}) such
that at $t=-1$ we have $h_i=B_iz$.  Then 
\begin{equation}\label{eq:limit.hi}
\lim_{i\to\infty} \sup_{\Sigma^i_j \cap B_2(0)} |\overline{b}_j z - h_i| = 0 \quad\text{for $j=1,2$,}
\end{equation}
 where $\overline{b}_j$ are the constants given in \eqref{eq:distinct}. 
 \end{proposition}
 
 \begin{remark}\label{rmk:Biz.approx.heat}
 The idea of Proposition \ref{prop:Biz.approx.heat} is that
 $\bar{b}_jz$ defines a solution of the heat equation on the union
 $P_1\cup P_2$, and we try to find solutions along the flows $L^i_t$
 which approximate it. Along the flows we do not have two components
 at each time converging to the two planes, so we cannot directly
 define a function like $\bar{b}_jz$. However,  in the
 limit as $i\to\infty$, the functions $B_i$ approximate the constants
 $\overline{b}_j$ on the two planes $P_j$. $B_i z$ only
 approximately satisfies the heat equation as $i\to\infty$ but
 should stay close to a genuine solution $h_i$ with the same initial
 condition. In addition we have a good pointwise estimate for the
 difference between $B_i$ and the constants $\overline{b}_j$ on the
 two components $\Sigma^i_j$ at the specific time $t=s_1$, as in Neves~\cite[Theorem
 4.2]{Neves:survey}. 
 \end{remark}

 \begin{proof}
 Let
\begin{equation}\label{eq:Ei}
 E_i = B_i z - h_i.
\end{equation}
Our goal is to show that $E_i$ is small as $i$ becomes large. At
$t=-1$ we have $E_i=0$, so we compute the evolution of $E_i$. We have
\[ \nabla (\beta_i + 2(t-s_1)\theta_i) = J( \bx^\perp + 2(s_1-t)\bH), \]
and since $\beta_i + 2(t-s_1)\theta_i$ satisfies the heat equation we
get
\[ (\partial_t - \Delta) B_i = |\bx^\perp + 2(s_1-t)\bH|^2 B_i. \]
Since $|B_i|\leq 1$,  at $t=-1$ we have $|h_i|\leq (1 +
|\mathbf{x}|^2)$. Using the maximum principle (see
Ecker--Huisken~\cite[Corollary 1.1]{EckerHuisken}, which applies to
subsolutions that satisfy the monotonicity formula) and the evolution
equation $(\partial_t - \Delta) (1 + |\mathbf{x}|^2) = -2n$ we find
that $|h_i| \leq C(1 + |\mathbf{x}|^2)$ for $t\in [-1,0)$ for a
dimensional constant $C > 0$. Below the constant $C$ may change from
line to line but is independent of $i, t$. In particular we also have
$|E_i| \leq C(1 + |\bx|^2)$.

Since $z$ and $h_i$ satisfy the heat equation along the flow $L^i_t$,
we have the evolution equation  
 \[
 (\partial_t -
   \Delta) E_i = |\bx^\perp + 2(s_1-t)\bH|^2 B_iz -2\langle\nabla B_i,\nabla z\rangle.
 \]
 We deduce that
   \[ |(\partial_t - \Delta) E_i| \leq |\bx^\perp +
     2(s_1-t)\bH|^2(1+|\bx|^2) + 2|\bx^\perp + 2(s_1-t)\bH|. \]
 From this, together with the estimate $|E_i| \leq C(1+|\bx|^2)$, we get
 \begin{equation}\label{eq:evol.Ei2}
   \begin{aligned}
     (\partial_t - \Delta) E_i^2 &\leq 2 |E_i| |(\partial_t -
     \Delta)E_i| - 2|\nabla E_i|^2 \\
     &\leq 2 |E_i| |\bx^\perp + 2(s_1-t) \bH|^2(1+|\bx|^2) + 4|E_i|\,
     |\bx^\perp + 2(s_1-t)\bH| \\
     &\leq E_i^2 + C (|\bx^\perp|^2 +
     (s_1-t)^2 |\bH|^2) (1 + |\bx|^4),
     \end{aligned}
   \end{equation}
   where we also used the estimate $4|E_i|b \leq E_i^2 + 4b^2$ to get the
   last line. 

   Using that $\theta_i$ satisfies the heat equation and
   $|\nabla\theta_i| = |\bH|$, as well as $(\partial_t - \Delta)
   |\bx|^4 \leq 0$, we also have
   \[ \begin{aligned}(\partial_t - \Delta) (1 + |\bx|^4)
     (t+1)\theta_i^2 &\leq (1+|\bx|^4)\theta_i^2 -
     2(1+|\bx|^4)(t+1)|\bH|^2  \\&\qquad - 4(t+1)\theta_i \langle \nabla
     \theta_i, \nabla |\bx|^4\rangle \\
     &\leq -(1+|\bx|^4) (t+1)|\bH|^2 + C(1 + |\bx|^4)\theta_i^2,
   \end{aligned} \]
   for $t\in (-1,0)$. 

   Let $\kappa > 0$ be small. Combining \eqref{eq:evol.Ei2} with
   the previous inequality, for $t\in (-1,0)$ we have
   \begin{equation}\label{eq:evol.Ei2.kappa}
     \begin{aligned} 
       (\partial_t - \Delta) \Big(e^{-t} E_i^2 &+ (1 +
       |\bx|^4)\kappa(t+1) \theta_i^2\Big) \\&
       \leq C(|\bx^\perp|^2+ (s_1-t)^2|\bH|^2) (1+|\bx|^4)  \\ &\quad+ \kappa
      C (1 + |\bx|^4)
       \theta_i^2 
       - \kappa(t+1) (1+|\bx|^4) |\bH|^2.
       \end{aligned} 
\end{equation}

   Suppose that $\bx_0 \in B_2(0) \cap L^i_{s_1}$ and denote by $\rho_{\bx_0, s_1}$
   the backwards heat kernel centred at $(\bx_0, s_1)$.
     For $t\in (-1,
   s_1)$ we have from  \eqref{eq:evol.Ei2.kappa}, using the
   monotonicity formula \eqref{eq:Huisken}, that
\begin{equation}\label{eq:evol.int.Ei2.kappa}
 \begin{aligned}
       \frac{d}{dt} \int_{L^i_t} (e^{-t}E_i^2 &+ (1 + |\bx|^4)\kappa
       (t+1)\theta_i^2) \rho_{\bx_0, s_1}\\ &\leq \int_{L^i_t} C (|\bx^\perp|^2 +
       (s_1-t)^2|\bH|^2) (1+|\bx|^4) \rho_{\bx_0, s_1} \\
       &\quad + \kappa C \int_{L^i_t} (1+|\bx|^4) \theta_i^2
       \rho_{\bx_0, s_1}\\
       &\quad - \kappa(t+1) \int_{L^i_t} |\bH|^2 (1+|\bx|^4)\rho_{\bx_0, s_1}.
     \end{aligned} 
     \end{equation}
    Integrating \eqref{eq:evol.int.Ei2.kappa} with respect to $t$ from
    $-1$ to $s_1$ yields: 
   \begin{equation}\label{eq:Ei1} 
   \begin{aligned}
    e^{-s_1}E_i^2(\bx_0,s_1) &+(1 + |\bx_0|^4) \kappa(s_1+1)\theta_i^2(\bx_0, s_1)
    \leq
     \int_{L^i_{-1}} e^{-t} E_i^2 \rho_{\bx_0, s_1} \\
     &+ \int_{-1}^{s_1} \int_{L^i_t} C (|\bx^\perp|^2 +
     (s_1-t)^2|\bH|^2) (1+|\bx|^4)\rho_{\bx_0, s_1}\rd t \\
     &+ \kappa C \int_{-1}^{s_1} \int_{L^i_t} (1+|\bx|^4)
     \theta_i^2\rho_{\bx_0, s_1} \\
     &-  \kappa\int_{-1}^{s_1} \int_{L^i_t} (t+1)|\bH|^2 (1 + |\bx|^4)\rho_{\bx_0, s_1} \rd t.
       \end{aligned}
     \end{equation}   
   Note that $E_i =0$ at $t=-1$ by the definition of $h_i$, and hence
   the first term on the right-hand side in \eqref{eq:Ei1} vanishes.
      
 We now estimate the second term in \eqref{eq:Ei1}.  Note that for $t\in [-1,s_1-\kappa]$ we have 
 \[
 (1 + |\bx|^4)\rho_{x_0,s_1}(\bx,t)\leq C_{\kappa}\rho_{0,0}(\bx,t)
 \]
 for some $\kappa$-dependent constant $C_{\kappa}>0$, since $s_1<0$ and thus $\rho_{\bx_0,s_1}$ will decay more rapidly at infinity than $\rho_{0,0}$ for any $t\in [-1,s_1-\kappa]$.  Therefore,
\begin{equation}\label{eq:second.term.est.1}
\begin{split}
 \int_{-1}^{s_1-\kappa}
     \int_{L^i_t}C (|\bx^\perp|^2 + 
     (s_1-t)^2|\bH|^2) &(1 + |\bx|^4)\rho_{\bx_0, s_1}\rd t\\
      &\leq C_\kappa \int_{-1}^{0} \int_{L^i_t} (|\bx^\perp|^2 +
     |\bH|^2)\rho_{0,0}\rd t.
     \end{split}
\end{equation} 
 We now notice that 
 \begin{equation}\label{eq:second.term.est.2} 
 \begin{split}
\int_{s_1-\kappa}^{s_1} \int_{L^i_t}C
     &(s_1-t)^2|\bH|^2 (1+|\bx|^4) \rho_{\bx_0, s_1}\rd t\\
     &\leq
                                                         C\kappa^2\int_{s_1-\kappa}^{s_1}
                                                         \int_{L^i_t}
                                                         |\bH|^2 (1+|\bx|^4)\rho_{\bx_0,s_1}\rd t\\
     &\leq 
       \frac{\kappa}{2}\int_{s_1-\kappa}^{s_1}\int_{L^i_t}(t+1)|\bH|^2(1
       + |\bx|^4)\rho_{\bx_0,s_1}\rd
       t,
       \end{split}
 \end{equation}
 where $C>0$ is a constant and $\kappa$ is chosen sufficiently small
 that $(t+1)\geq 2C\kappa$ for $t\in [s_1-\kappa,s_1]$.  Equation
 \eqref{eq:second.term.est.2} shows that the integral on the left-hand
 side of the inequality can be compensated for using the last term in
 \eqref{eq:Ei1}.

 Our remaining concern is 
 \begin{equation}\label{eq:second.term.est.3}
 \begin{aligned}
  \int_{s_1-\kappa}^{s_1} \int_{L^i_t}& |\bx^\perp|^2  (1+|\bx|^4)
  \rho_{\bx_0, s_1}\rd t \leq
  \int_{s_1-\kappa}^{s_1}\int_{L^i_t}|\bx|^2  (1+|\bx|^4)  \rho_{\bx_0,s_1}\rd t\\
  &=\int_{s_1-\kappa}^{s_1}\int_{L^i_t\cap
    B_{\kappa^{-1/10}}(0)}|\bx|^2  (1+|\bx|^4) \rho_{\bx_0,s_1}\rd t\\
    &\quad +\int_{s_1-\kappa}^{s_1}\int_{L^i_t\setminus
    B_{\kappa^{-1/10}}(0)}|\bx|^2  (1+|\bx|^4) \rho_{\bx_0,s_1}\rd t. 
  \end{aligned}
 \end{equation}
The first integral can clearly be estimated as
\begin{equation}\label{eq:second.term.est.4}
\begin{split}
 \int_{s_1-\kappa}^{s_1}\int_{L^i_t\cap B_{\kappa^{-1/10}}(0)}&|\bx|^2
(1+|\bx|^4) \rho_{\bx_0,s_1}\rd t
\\ &\leq
2\kappa^{-6/10}\int_{s_1-\kappa}^{s_1}\int_{L^i_t\cap
  B_{\kappa^{-1/10}}(0)}\rho_{\bx_0,s_1}\rd t\leq C\kappa^{2/5} 
  \end{split}
\end{equation}  
for some constant $C>0$, using the uniform area bounds for $L^i_t$.
Using the area bounds again for
$t\in [-1,s_1]$, we can estimate our remaining spacetime
integral by the integral over an $n$-plane $P$ for $\kappa$ sufficiently
small: 
\begin{equation}\label{eq:second.term.est.5}
\begin{split}
  \int_{s_1-\kappa}^{s_1}\int_{L^i_t\setminus
    B_{\kappa^{-1/10}}(0)}&|\bx|^2  (1+|\bx|^4) \rho_{\bx_0,s_1}\rd t\\
    &\leq
  C_1 \int_{-\kappa}^0\int_{P\setminus B_{\kappa^{-1/10}(0)}}|\bx|^2
  (1+|\bx|^4) \rho_{\bx_0,0}\rd t\\
  &\leq C_2 e^{-1/\kappa} 
  \end{split}
\end{equation}
for constants $C_1,C_2>0$. 

Combining \eqref{eq:second.term.est.1}--\eqref{eq:second.term.est.5} shows that, for $\kappa$ sufficiently small, we have  
\begin{equation}\label{eq:second.term.est}
 \begin{aligned} \int_{-1}^{s_1} \int_{L^i_t}C &(|\bx^\perp|^2 +
   (s_1-t)^2|\bH|^2)  (1+|\bx|^4) \rho_{\bx_0, s_1}\rd t\\
   &
   \leq  
     C_\kappa \int_{-1}^{0} \int_{L^i_t} (|\bx^\perp|^2 +
     |\bH|^2)\rho_{0,0}\rd t \\
  &\quad +C\kappa^{2/5}+\frac{\kappa}{2} \int_{s_1-\kappa}^{s_1}\int_{L^i_t}(t+1)|\bH|^2  (1+|\bx|^4) \rho_{\bx_0,s_1}\rd t
   \end{aligned}
   \end{equation}
   for some constant $C>0$ and a constant $C_{\kappa}>0$ depending on $\kappa$.   
 
 Noting also that $\theta_i^2$ is uniformly bounded, we may therefore combine \eqref{eq:Ei1} and \eqref{eq:second.term.est} to obtain 
\begin{equation}\label{eq:Ei2}
 E_i^2(\bx_0, s_1) \leq C_\kappa\int_{-1}^0 \int_{L^i_t} (|\bx^\perp|^2
   + |\bH|^2)\rho_{0,0}\rd t + C\kappa^{2/5},
   \end{equation}
 if $\kappa > 0$ is sufficiently small.
 
 For fixed $\kappa >
 0$, the first term on the right-hand side of \eqref{eq:Ei2} converges
 to zero as $i\to \infty$, as in \cite[Lemma 5.4]{Neves.ZM}. It 
 follows that for any $\kappa > 0$ we can choose $i$ sufficiently
 large so that 
 $$E_i^2(\bx_0, s_1) \leq 2C\kappa^{2/5}.$$ By definition of $E_i$ in
 \eqref{eq:Ei}, and the fact that $B_i=\cos(\beta_i)$ at $t=s_1$, we 
 have that
\begin{equation}\label{eq:lim.Ei}
\lim_{i\to\infty} \sup_{B_2(0) \cap L^i_{s_1}} |\cos(\beta_i) z-h_i| = 0. 
\end{equation}

As in \cite[Lemma 7.3]{Neves.ZM}, we now use that the limiting
behaviour of the functions
$B_i$ in \eqref{eq:Bi} as $i\to\infty$ is $t$-independent. More
precisely, for all $\phi$ with compact support in $B_2(0)$  and $f\in
C^2(\mathbb{R})$ we have
\begin{equation}\label{eq:Bistatic}
  \lim_{i\to\infty} \int_{L^i_{s_1}} f(B_i) \phi\, d\mathcal{H}^n =
  \lim_{i\to\infty} \int_{L^i_{-1}} f(B_i) \phi\, d\mathcal{H}^n.
\end{equation}
On $L^i_{s_1}$ we have $B_i = \cos(\beta_i)$  and so we have the
pointwise bound $|\nabla B_i|\leq |\bx^\perp|$. Using the Poincar\'e
type inequality  \cite[Proposition A.1]{Neves.ZM}, we deduce that
there are constants $\hat{b}_1, \hat{b}_2$ such that
$\sup_{\Sigma^i_j\cap B_2(0)} |B_i - \hat{b}_j| \to 0$ as $i\to\infty$. At the
same time from \eqref{eq:Bistatic} we find that $\hat{b}_j =
\overline{b}_j$ for the constants in \eqref{eq:distinct}, since on
$L^i_{-1}$ we have $B_i = \cos(\beta_i - 2(1+s_1)\theta_i)$. Note that
by construction on $L^i_{-1}$ we have 
$\beta_i\to \overline{\beta}_j$ and $\theta_i \to \overline{\theta}_j$
on the plane $P_j$, locally smoothly away from $\ell$. It follows
then from \eqref{eq:lim.Ei} that $\lim_{i\to\infty}
\sup_{\Sigma^i_j\cap B_2(0)}
|\overline{b}_j z - h_i| = 0$, as required. 
 \end{proof}

\subsection{The linking argument} \label{sec:linking-argument} Continuing the setup from the
previous subsection, 
we  now show that indeed the second possibility in
Proposition~\ref{prop:2cases} leads to a contradiction, if our flow
$t\mapsto L_t$ is smooth and embedded. 

\begin{proposition}
\label{prop:linking}
Suppose that we have $\phi_i\in \Span\{x_1,\ldots,x_{2n-2},z\}$ with $\phi_i\to 0$ 
 and a sequence
$\lambda_i\to 0$ such that along the sequence $L^i_t$ we have
\begin{equation}\label{eq:ui.a.limit}
u_i=\lambda_i^{-1}(w-\phi_i)\to z\theta\quad\text{on  $P_1\cup P_2$ as $i\to\infty$,}
\end{equation}
where the convergence is in $L^2$ and locally uniform away from
$\ell$. Then for sufficiently large $i$ the flow $L^i_t$ is not embedded. 
\end{proposition}

\begin{remark}  Recall that $\Sigma^i_j$ are the components of
  $B_2(0)\cap L^i_{s_1}$, as in Definition \ref{dfn:scale.cpts}, and
  let us suppose for simplicity that $C^i_j=\Sigma^i_j\cap \partial
  B_1(0)$ are smooth $(n-1)$-dimensional submanifolds of the sphere
  (which can always be done by changing the
  radius of the ball slightly if necessary). 
The key to the argument is to 
show that the submanifolds $C^i_j$ in the $(2n-1)$-sphere are 
linked for $i$ sufficiently large, which implies that the 
$\Sigma^i_j$ intersect in $B_2(0)$. Then $L^i_{s_1}$ cannot be
embedded.  
\end{remark}

\begin{proof}
Since $\theta$ equals the distinct constants $\overline{\theta}_1,
\overline{\theta}_2$ on the planes $P_1,P_2$, by modifying the
$\lambda_i$ and adding multiples of $z$ to the $\phi_i$,
we can assume for simplicity that 
\[
u_i\to \overline{b}_jz\quad\text{on $P_j$ as $i\to\infty$,}
\]
where $\overline{b}_j$ are given in \eqref{eq:distinct}. The
convergence is in $L^2$, and locally uniform away from the line
$\ell$.  We  also
assume without loss of generality that $\lambda_i > 0$. 

Recall the notation of Definition~\ref{dfn:scale.cpts} and
Proposition~\ref{prop:Biz.approx.heat}.  At $t=-1$ we have $h_i =
B_iz$ by definition, and the function $B_i$ converges in $L^2$ and
locally smoothly away from $\ell$ to the constants $\bar{b}_j$ on the
two planes. It follows that at $t=-1$ we have $\Vert u_i -
h_i\Vert_{L^2} = 0$. The
monotonicity formula applied with points $(\bx_0, s_1)$ for different
$\bx_0\in B_2(0)$ then implies 
   \[ \lim_{i\to\infty} \sup_{\Sigma^i_j} |u_i - h_i| =
     0. \]
 Applying Proposition \ref{prop:Biz.approx.heat} then yields
   \[  \lim_{i\to\infty} \sup_{\Sigma^i_j} |u_i - \overline{b}_j z| =
     0. \]
 We deduce that, given any $\epsilon > 0$, once $i$ is sufficiently
   large we will have
\begin{equation}\label{eq:w.eps.ineq} |w - \phi_i - \lambda_i \overline{b}_j z| <
     \epsilon \lambda_i \quad\text{on $\Sigma^i_j$.}
     \end{equation}

   Suppose without loss of generality that $\overline{b}_1>\overline{b}_2$  and choose $$0<\epsilon < |\overline{b}_1-\overline{b}_2|/100$$ in \eqref{eq:w.eps.ineq}.  Let
   $$\overline{b}_0 = (\overline{b}_1 + \overline{b}_2)/2$$
   and, recalling that $\phi_i  \in\Span\{x_1,\ldots,x_{2n-2},z\}$, define half-spaces
\[\begin{aligned} 
 \mathcal{H}^i_+&=\{(x_1,\ldots,x_{2n-2},z,w)\in\CC^n\,:\,w >
   \phi_i+ \lambda_i\overline{b}_0 z\},
   \\
   \mathcal{H}^i_-&=\{(x_1,\ldots,x_{2n-2},z,w)\in\CC^n\,:\,w <
   \phi_i+ \lambda_i\overline{b}_0z\}.
   \end{aligned}
   \]    
   The inequality \eqref{eq:w.eps.ineq} implies that, for all $i$ sufficiently large,
 \begin{equation}\label{eq:halfspaces}
 \begin{gathered}
 (\Sigma^i_1\cap\{z>1/2\})\subseteq\mathcal{H}^i_+\quad\text{and}\quad(\Sigma^i_2\cap\{z>1/2\})\subseteq\mathcal{H}^i_-,\\
 (\Sigma^i_1\cap\{z<-1/2\})\subseteq\mathcal{H}^i_-\quad\text{and}\quad(\Sigma^i_2\cap\{z<-1/2\})\subseteq\mathcal{H}^i_+.
 \end{gathered}
 \end{equation}
 In other words, the relative positions of the components $\Sigma^i_j$ in terms of the halfspaces $\mathcal{H}^i_{\pm}$ must switch as we pass from $z>1/2$ to $z<-1/2$.

We can choose $R=1+\delta$ for $\delta\geq 0$ small such that 
\begin{equation}\label{eq:curves}
C_j^i = \Sigma^i_j \cap \partial B_R(0)
\end{equation}
 are smooth.  Our aim now is to show that the
 submanifolds  $C_j^i$ in  $\partial B_R(0)$ are linked for
 sufficiently large $i$, which will imply that the $\Sigma^i_j$
 intersect in $B_R(0)$.
 
   Consider the
   two points $p_-, p_+$ whose only non-zero entries are $\pm R$ in
   the $z$-component in coordinates
   $(x_1,\ldots,x_{2n-2},z,w)$ on $\CC^n$. So $p_-, p_+$  lie on $\ell\cap \partial
   B_R(0)$ where $\ell=P_1\cap P_2$.
   Since the $\Sigma^i_j$ converge smoothly to $P_j$ away from the singular
   line $\ell$ as $i\to\infty$, we can assume that outside of $B_{1/20}(p_\pm)$ the submanifolds
   $C_j^i$ are smooth perturbations of $P_j\cap \partial B_R(0)$. 
   
Any connected components of the $C_j^i$ contained entirely
inside $B_{1/10}(p_{\pm})$ must lie in different half-spaces
$\mathcal{H}^i_{\pm}$ by \eqref{eq:halfspaces} for $i$ sufficiently
large, and so do not contribute to the linking number of the 
$C_j^i$.  We may therefore discard these components, if there are any,
and assume from now on that the $C^i_j$ are connected.    
   
For $j=1,2$, let $\tilde{P}^i_j$ be the graph of $w = \phi_i +
   \lambda_i \overline{b}_j z$ over
   $P_j$. Since $\phi_i, \lambda_i\to 0$, the 
   $\tilde{P}^i_j$ are small  perturbations of the $P_j$ for $i$
   sufficiently large.  Moreover, since
   $\overline{b}_1\neq\overline{b}_2$ by \eqref{eq:distinct} we have
   that the $\tilde{P}^i_j$ intersect transversely at the origin.
   Hence, the spheres
\begin{equation}\label{eq:tilde.curves}   
   \tilde{C}^i_j = \tilde{P}^i_j\cap \partial
   B_R(0)
   \end{equation} have linking number 1.
   
We now claim that, for $i$ sufficiently large, the submanifolds $C^i_j$ in
\eqref{eq:curves} can be deformed to the $\tilde{C}^i_j$ in
\eqref{eq:tilde.curves} without any crossings.  Outside of 
   the balls $B_{1/20}(p_\pm)$ this is clear since there the $C^i_j$ 
   are smooth perturbations of the $\tilde{C}^i_j$.  At the same time,
   for $i$ sufficiently large, 
   inside the balls $B_{1/10}(p_\pm)$ the pairs of submanifolds $\{ C^i_j, 
   \tilde{C}^i_j\}$ are contained in disjoint half-spaces for $j=1,2$
   by \eqref{eq:halfspaces}, so the submanifolds in each pair can be
   deformed to 
   coincide without intersecting the submanifolds in the other pair.  
   
   We conclude that the linking
   number of the submanifolds $C^i_j$ in \eqref{eq:curves} is therefore
   also 1 for sufficiently large $i$, which implies that the flow is
   not embedded. 
\end{proof}

\section{Proof of Main Theorem}

We first show that combining the results from Section \ref{sec:three-annulus} and Section \ref{sec:linking-argument} yields a proof of the main theorem:

\begin{proof}[Proof of Theorem \ref{thm:translators.intro}] We can assume that the ancient flow $\cM$ is defined for $t<0$. We note that the assumption that the flow has a blow-down given by the static union of the planes $P_1\cup P_2$ implies that the entropy is bounded above by 2. This implies that if the flow has an immersed point $(\bx_0,t_0)$, then the monotonicity formula yields that the flow is backwards self-similar around $(\bx_0,t_0)$, i.e.~the flow is given by the static flow $(\cM_{P_1\cup P_2} + (\bx_0,t_0))\cap\{t<0\}$.

We can thus assume that $\cM$ is embedded. Combining Proposition~\ref{prop:2cases} and Proposition~\ref{prop:linking} yields the statement.
\end{proof}

Since in many geometric applications it is essential to classify not only smooth ancient solutions to mean curvature flow, but also more general Brakke flows arising as limit flows, we also record the following extension of Theorem \ref{thm:translators.intro}.

\begin{theorem} \label{thm:translators-extension} 
 Let $P_1,P_2\subset \CC^2$ be Lagrangian subspaces which intersect along a line $\ell$ through $0$.  Let $\cM$ be an ancient $2$-dimensional Brakke flow which is the (weak) limit of smooth, zero-Maslov, exact Lagrangian mean curvature flows $(L^i_t)_{-R_i^2<t<0}$ defined on $B(0,R_i) \subset \CC^n$, where $R_i \to \infty$, with uniformly bounded variation of the Lagrangian angle and uniformly bounded area ratios. If $\cM$ has a blow-down at $-\infty$ given by the static flow consisting of the union of the planes $P_1\cup P_2$, then $\cM$ is a smooth translator.
\end{theorem}

\begin{proof} We again assume that $\cM$ is defined and non-vanishing for $t<0$. Note that the assumptions imply that $\cM$ has uniformly bounded area ratios and is unit regular, meaning that every point of Gaussian density one has a space-time neighbourhood where the flow is smooth. 
  Furthermore, as in the proof of Theorem \ref{thm:translators.intro}, it follows that the entropy is bounded above by 2. Assume now that there is a point $(\bx_0,t_0)$ where the Gaussian density of $\cM$ is 2. Then as above we see that $\cM=(\cM_{P_1\cup P_2} + (\bx_0,t_0))\cap\{t<0\}$ (using unit regularity to conclude that neither of the two planes can vanish before $t=0$).

We can thus assume that all Gaussian density ratios of $\cM$ are strictly less than 2.   Neves structure theory \cite{Neves.ZM} then implies that we obtain uniform local curvature bounds along the sequence $(L^i_t)_{-R_i^2<t<0}$, and thus the convergence is smooth. This yields that $\cM$ is a smooth, ancient, zero-Maslov, exact Lagrangian mean curvature flow with uniformly bounded variation of the Lagrangian angle and uniformly bounded area ratios.   Theorem \ref{thm:translators.intro} then implies the statement.
\end{proof}

\begin{remark}
 In the previous theorem one can also allow that the flows $(L^i_t)_{-R_i^2<t<0}$ are defined on the Riemannian manifolds $(B(p_i,R_i), g_i)$ where $B(p_i,R_i) \subset \RR^{2n}$ are geodesic balls with respect to $g_i$ and $g_i$ is a sequence of Calabi--Yau metrics on $B(p,R_i)$ converging smoothly to the standard Euclidean metric on $\CC^n$.
\end{remark} 

\begin{appendix}
\section{Smoothly immersed Brakke flows}
\label{app:monotonicity}

\noindent In this appendix we give a definition of smoothly immersed Brakke flows and record that the weighted version of Huisken's monotonicity formula, see \cite[Theorem 4.13]{Ecker.book}, can be extended to weights with polynomial growth. We recall Remark \ref{rem:smoothly_immersed} for the relation with properly immersed mean curvature flows. 

\begin{definition}[Smoothly immersed Brakke flows]\label{def:smoothly_immersed} We say that an $n$-dimensional (integral) Brakke flow $\cM$ in $\RR^{n+m}$, defined on a time interval $I\subset \RR$, is \emph{smoothly immersed} if around every point $(\bx, t) \in \RR^{n+m}\times I$ in the support $\supp\cM$ of the Brakke flow there exists an open space-time neighbourhood $U$, such that the flow can be represented by a smoothly immersed mean curvature flow in $U$. Note that this includes that the multiplicity $\theta(\bx, t)$ agrees with the number of sheets passing through $(\bx, t)$. If $n=m$ we say that $\cM$ is in addition \emph{Lagrangian} if the local immersions can be chosen to be Lagrangian.
\end{definition}

Consider an $n$-dimensional, smoothly immersed Brakke flow $\cM$ in $\RR^{n+m}$. We call a map $f$ a \emph{function} on $\cM$ if it assigns to each $(\bx, t) \in \supp\cM$ an unordered $\theta(\bx,t)$-tuple of real numbers. We say that such a function is \emph{continuous, smooth, etc.}~if locally it can be represented by a \emph{continuous, smooth, etc.}~function on a suitable local immersion representing $\cM$. 

\begin{definition}[Functions with polynomial growth]\label{def:polynomial_growth} We say that such a function $f$ on $\cM$ has \emph{polynomial growth} $d \in \NN$ if for every interval $J\Subset I$ there exists a constant $C_J$ such that for all $R>0$
\begin{equation}\label{eq:polynomial_growth}
\sup_{\supp\cM\cap (B_R(\bOh)\times J)} \|f(\bx, t)\| \leq C_J (1+R^d)\, ,
\end{equation}
where $\|f(\bx, t)\|$ denotes the maximal possible absolute value of $f$ at $(\bx,t)$.
\end{definition}

We recall that Ecker's weighted version of Huisken's monotonicity formula is given by \eqref{eq:Huisken} for a (sufficiently smooth) function $f$ on $\cM$ with (uniform) compact support.  We can extend this as follows to a (sufficiently smooth) function with non-compact support. Note that in the following we will do all calculations implicitly on the local immersions representing the flow.

\begin{proposition}\label{prop:monotonicity} Let $\cM$ be an $n$-dimensional, smoothly immersed Brakke flow $\cM$ in $\RR^{n+m}$ and let $f$ be a smooth function on $\cM$. Let $t_1, t_2 \in I$ with $t_1<t_2<t_0$ and assume that 
 $$ \int_{t_1}^{t_2} \int \left(f^2 + \left|\left(\partial_t -\Delta\right) f\right|\right) \rho_{\bx_0,t_0}\, d\mu_t\, dt < \infty\, ,$$
 as well as $f(\cdot, t_1) \in L^1(\rho_{\bx_0,t_0}(\cdot, t_1)\, d\mu_{t_1})$ and $f(\cdot, t_2) \in L^1(\rho_{\bx_0,t_0}(\cdot, t_2)\, d\mu_{t_2})$. Then
 \begin{equation}\begin{split}
  \int f\, \rho_{\bx_0,t_0}\, d\mu_{t_2} &\leq \int f\, \rho_{\bx_0,t_0}\, d\mu_{t_1} + \int_{t_1}^{t_2} \int \left(\partial_t -\Delta\right) f\, \rho_{\bx_0,t_0}\, d\mu_t\, dt\\ &\qquad - \int_{t_1}^{t_2} \int f \left| \bH - \frac{(\bx-\bx_0)^\perp}{2(t-t_0)}\right|^2\, \rho_{\bx_0,t_0}\, d\mu_t\, dt\, .
  \end{split}
  \end{equation}
 \end{proposition}

\begin{proof}
 Let $\eta$ be an ambient cut-off function. We have
 $$ \left(\partial_t-\Delta\right) (\eta f) = \eta \left(\partial_t -\Delta\right) f + f \left(\partial_t -\Delta\right) \eta - 2\langle \nabla f, \nabla \eta\rangle $$
 and thus after integration by parts
 \begin{equation}\label{eq:monotonicity.1}
\begin{split}
\int \left(\partial_t -\Delta\right) (\eta f)\, \rho \,d\mu_t &= \int \eta \left(\partial_t  -\Delta\right) f  \rho \, d\mu_t \\
&\qquad +  \int  f\left( \left(\partial_t  +\Delta\right) \eta
  + 2 \langle \nabla \eta, \tfrac{\nabla \rho}{\rho}\rangle\right)\rho \, d\mu_t\, .
 \end{split}
\end{equation}
Recall that for an ambient function $\eta$  one has at $(\bx, t) \in \supp \cM$
\begin{equation}\label{eq:ambient_laplace}
 \Delta_{M^i_t} \eta = \text{tr}_{T_pM^i_t} D^2\eta + \langle D\eta, \bH_{M^i_t}\rangle\, ,
\end{equation}
where $D$ is the standard ambient derivative and $M^i_t$ is one of the sheets locally representing the flow. This implies (assuming $\eta$ is independent of time)
$$ \left(\partial_t  +\Delta_{M^i_t}\right) \eta = \text{tr}_{T_pM^i_t} D^2\eta + 2\langle D\eta, \bH_{M^i_t}\rangle\, .$$
Recall further that the assumption of bounded area ratios implies that
$$ \int_{t_1}^{t_2} \int \left| \bH - \frac{\mathbf{x}^\perp}{2t}\right|^2\, \rho\, d\mu_t\, dt \leq C<\infty,$$
where $\bH(\bx,t) = \sum_{i=1}^{\theta(\bx,t)} \bH_{M^i_t}(\bx,t)$ is the varifold mean curvature and the $M^i_t$ are the sheets passing through $(\bx,t)$. Thus, again using bounded area ratios, we have 
$$ \int_{t_1}^{t_2} \int\left| \bH\right|^2\, \rho\, d\mu_t\, dt \leq C(t_2)\, .$$
The above gives
\begin{equation}\label{eq:monotonicity.2}
\begin{split}
 \bigg|\int  f&\left( \left(\partial_t  +\Delta\right) \eta + 2 \langle \nabla \eta, \tfrac{\nabla \rho}{\rho}\rangle\right)\rho \, d\mu_t\bigg|\\ 
 &\leq C \int |f|  \left(|D^2\eta| + |\nabla \eta| \tfrac{|\nabla \rho|}{\rho} + |D\eta| |\bH|\right) \rho\, d\mu_t\\
 &\leq C \left(\int |f|^2 \rho\, d\mu_t\right)^\frac{1}{2} \left(\int  \left(|D^2\eta|^2 + |\nabla \eta|^2 \tfrac{|\nabla \rho|^2}{\rho^2} + |D\eta|^2 |\bH|^2\right) \rho\, d\mu_t\right)^\frac{1}{2}
\end{split}
 \end{equation}
 We now choose $\varphi$ to be a cutoff function which is equal to one on $B_1(0)$ and vanishes outside of $B_2(0)$ and let $\eta_R(\mathbf{x}) = \varphi(\mathbf{x}/R)$. Integrating \eqref{eq:monotonicity.1} from $t_1$ to $t_2$ with $\eta = \eta_R$ and letting $R\to \infty$ (using \eqref{eq:monotonicity.2} and the space-time integral bound on $|\bH|^2$) gives the result.
\end{proof}

\begin{remark} \label{rem:mon_poly} For a  smooth function $f$ on $\cM$ with polynomial growth such that $(\partial_t  - \Delta)f$ also has polynomial growth the conditions of Proposition \ref{prop:monotonicity} are satisfied.
\end{remark}

We note that using polynomial barriers one obtains existence of caloric functions with polynomial growth given initial data of polynomial growth. Furthermore, the above monotonicity formula yields uniqueness. We record this in the following proposition. We write $\cM_{t\geq t_0}$ for the restriction of the Brakke flow to $I\cap \{t\geq t_0\}$ and extend the definition of a function with polynomial growth on $\supp (\cM)\cap \{(\bx, t_0)\, |\, \bx \in \RR^{n+m}\}$ in the obvious way.

\begin{proposition} \label{prop:caloric_unique} Let $\cM$ be an $n$-dimensional, smoothly immersed Brakke flow  in $\RR^{n+m}$ and for $t_0\in I$ let $f_0$ be a smooth function on $\supp (\cM)\cap \{(\bx, t_0)\, |\, \bx \in \RR^{n+m}\}$ with polynomial growth. Then there exists a unique smooth function $f$ on $\cM_{t\geq t_0}$ of polynomial growth such that 
\begin{equation}\label{eq:heat_equation}
(\partial_t -\Delta) f = 0 \quad \text{and}\quad  f|_{t=t_0} = f_0\, .
\end{equation}

 \end{proposition}
 \begin{proof}
We first note that Proposition \ref{prop:monotonicity} together with Remark \ref{rem:mon_poly} yields uniqueness as stated. For existence we have the following claim.\\[1ex]
{\bf Claim:} Given $R>0$ there exists a solution $f_R$ to \eqref{eq:heat_equation} on $(B_R(\bOh)\times [t_0, t_0+R^2)) \cap \supp \cM$ such that for every $0<r\leq R$
$$\sup_{\supp \cM\cap (B_{r}(\bOh)\times [t_0,t_1])} \|f_R\| \leq C_0e^{C_1(t_1-t_0)} (1+r^d)\, ,$$
for some $C_1>0$ just depending on $n,m$ and $d$.\\[1ex]
Note that the claim does not specify any boundary values at the spatial boundary. From the claim, the existence follows, since interior higher order estimates imply that we can take a subsequential limit as $R \to \infty$ to obtain the stated solution $f$. The convergence of $f$ to the initial data $f_0$ as $t\searrow t_0$ follows similarly from higher order interior estimates.

To prove the claim, let $R>0$ be given. Note that since $\cM$ is smooth there exists $K>0$ such that the mean curvature $\bH$ of the flow is bounded by $K>0$ on $B_{4R}(\bOh)\times [t_0, t_0+R^2]$. Since $\cM$ is smoothly immersed, there exists an (open) $n$-manifold $M$ and an immersion $F_{t_0}$ such that $F_{t_0}:M\to \RR^{n+m}$ smoothly parametrises $\cM(t_0)\cap B_{3R}(\bOh)$. Furthermore we can smoothly extend $F_{t_0}$ to a standard (immersed) mean curvature flow $F_t$ parametrising $\cM$ for $t\in [t_0, t_0 + \delta]$ with $\delta:= R/(2K)$: this follows from the bound on $\bH$. Note further that $\cM(t+\delta)\cap B_R(\bOh) \subset F_{t+\delta}(M)$. This again follows from the bound on $\bH$ on $B_{4R}(\bOh)\times [t_0, t_0+R^2]$. 

Let $f_0$ be of polynomial growth of degree $d$ such that for all $r>0$
\begin{equation}\label{eq:inital_data}
 \sup_{(\bx, t_0) \in \supp\cM\cap (B_r(\bOh)\times \{t_0\})} \|f_0(\bx)\|  \leq C_0 (1+r^d)\, .
\end{equation}
Choose $R'\in (2R, 3R)$ such that $U:=F_{t_0}^{-1}(B_{R'}(\bOh)) \subset M$ has smooth boundary. We can then construct a solution $\hat f$ to the heat equation (with respect to the induced metric $g_t$ on $M$ via $F_t$) with initial value $f_0$ and boundary value zero. Note that \eqref{eq:ambient_laplace} together with $d\geq 2$ implies that there exists $C_1>0$ such that $C_0e^{C_1(t-t_0)}(1+|\bx|^d)$ is a supersolution to the heat equation along the flow.  Thus,  by the maximum principle together with \eqref{eq:inital_data}, as well as the assumptions on the boundary data, we have 
$$ |\hat f(x,t)| \leq C_0e^{C_1(t-t_0)} (1+|F(x,t)|^d)\, ,$$
for $(x,t) \in U\times [t_0,t_0+\delta]$. Note that by restriction this yields a solution $f$ to the heat equation with the claimed bounds on $\cM \cap (B_R(\bOh)\times [t_0,t_0+\delta])$. We can now repeat this process, starting at $t_0+\delta$ where as initial data we take $\hat f \circ F^{-1}_{t_0+\delta}$ on $F_{t_0+\delta}(U) \subset \supp \cM(t_0+\delta)\cap B_{4R}(\bOh)$ and zero on $(\supp \cM(t_0+\delta)\cap B_{4R}(\bOh)) \setminus F_{t_0+\delta}(U)$. Repeating this process finitely many times yields the stated solution $f_R$.
\end{proof}

\end{appendix}


\providecommand{\bysame}{\leavevmode\hbox to3em{\hrulefill}\thinspace}
\providecommand{\MR}{\relax\ifhmode\unskip\space\fi MR }
\providecommand{\MRhref}[2]{%
  \href{http://www.ams.org/mathscinet-getitem?mr=#1}{#2}
}
\providecommand{\href}[2]{#2}

\end{document}